\documentclass[a4paper,11pt,twoside]{article}
\usepackage{pst-fill,pst-grad}
\usepackage{textcomp}
\usepackage[english]{babel}
\usepackage[utf8]{inputenc} 
\usepackage[titletoc]{appendix}
\usepackage{titlesec}
\usepackage{graphicx}
\usepackage{amsmath}
\usepackage{float} 
\usepackage{fancyhdr} 
\usepackage[matrix,arrow,curve]{xy} 
\usepackage{pstricks} 
\usepackage{amsmath,amsfonts,verbatim,afterpage,theorem,euscript,mathrsfs,amssymb}
\usepackage{amsfonts}
\usepackage{amssymb}
\usepackage{array}
\usepackage{dsfont}
\usepackage[colorlinks=true,linkcolor=blue,citecolor=red]{hyperref}
\usepackage{authblk}
\usepackage{color}
 
\newtheorem{Proposition}{Proposition}[section]

\newtheorem{Lemme}{Lemma}[section]

\newtheorem{TheoremeP}{Theorem}

\newtheorem{Remarque}{Remark}
\newcommand{\sign}{\text{sign}}

\def \R{\mathbb{R}}
\def \Rn{\mathbb{R}^n}
\def \Rt{\mathbb{R}^3}

\def \finpv{\hfill $\blacksquare$  \\ \newline }
\def \pv{{\bf{Proof.}}~} 

\def \ds{\displaystyle}
\def \rd{\textcolor{black}}
\def \bl{\textcolor{black}}  

\title{\bf Spatial behavior of solutions for a large class of non-local PDE's arising from stratified flows.}

\author[1]{ Manuel Fernando Cortez\footnote{ manuel.cortez@epn.edu.ec}}
\author[2]{ Oscar Jarr\'in\footnote{corresponding author: oscar.jarrin@udla.edu.ec}}

\affil[1]{\scriptsize Departamento de Matem\'aticas, Escuela Politécnica Nacional,  Ladron de Guevera E11-253, Quito, Ecuador}
\affil[2]{\scriptsize Dirección General  de Investigación  (DGI),
	Universidad de las Américas,
	Calle José Queri s/n y Av. Granados. Bloque 7, Tercer Piso, Quito, Ecuador.} 
\date{\today}
\begin{document}
\maketitle	
\begin{abstract}
We propose a  theoretical model of a non-local dipersive-dissipative equation which contains as a particular case a large class of non-local PDE's arising from stratified flows. Within this fairly general framework, we study the spatial behavior of solutions proving some sharp pointwise and averaged decay properties as well as some pointwise grow properties. \\[3mm]
\end{abstract}
\section{Introduction} 
Stratified flows, which roughly speaking are fluids with a density variation, are everywhere in nature and play a key role in a range of natural phenomena, from ocean circulation to weather forecasting. Mathematical models for these flows  are helpful in understanding the real world. These models essentially write down as non-local,  dispersive-dissipative type equations, see \emph{e.g.} the range of  equations (\ref{OST})-(\ref{dgBO}) below, and these equations  describe the evolution of  nonlinear internal long waves considering  different physical settings. We refer to  \cite{Abde,Benjamin,Benjamin1,Benjamin2,OST0,OST,OST1,Saut} and the references therein for a small sample of the huge existing literature.\\

In this article, we propose a \emph{theoretical} equation which contains as a particular case some well-known relevant physical model arising from stratified fluids. Within the fairly general setting of this  equation, we investigate some \emph{sharp properties} of the spatial behavior of solutions. \\
 
Let us consider the following  Cauchy problem for a dispersive-dissipative equation with a non-local perturbation term:   
\begin{equation}\label{Equation} 
\left\{\begin{array}{ll}\vspace{2mm} 
\partial_t u + D (\partial_x u) +u^k \partial_x u + \eta(\mathcal{H} \partial^{n}_x u +   \mathcal{H}_m u)  = 0, & \eta >0,  \\
u(0, \cdot) = u_0. & 
\end{array} 
\right.
\end{equation}
In this equation,  the \emph{dispersion} effects are given by the term  $\ds{D (\partial_x u)}$, where  $D$  is a  pseudo-differential operator $D$ defined in the Fourier variable as follows: for $\varphi \in \mathcal{S}(\R)$ 
\begin{equation}\label{def-P}
\widehat{D (\varphi)}(\xi)= p(\xi) \widehat{\varphi}(\xi).
\end{equation}   The symbol $p(\xi)$ characterizes the linearized dispersion relation of the model equation (\ref{Equation}). We consider here a fairly general symbol $p(\xi)$ verifying the following natural assumptions  (see the Section $5$ of \cite{Abde}):  $p\in L^{\infty}_{loc}(\R)$  is a \emph{real-valued} function, continuous at the origin, smooth outside the origin  and with polynomial grow at infinity, \emph{i.e.},  for \emph{a.e.} $\xi \in \R$ we have 
\begin{equation}\label{estim-symbol}
\vert p(\xi) \vert \leq c \vert \xi \vert^{\sigma}, \quad \text{with}\quad \sigma >0. 
\end{equation}
It is easy to observe that the operator $D$ commutes with differentiation and moreover, since the symbol $p(\xi)$ is real-valued then $D$ is a self-adjoint operator on its domain in the space $L^2(\R)$.\\

Thereafter, for a parameter  $k\in \mathbb{N}^{*}$ we consider a fully non-linear term of the form $u^k \partial_x u$. Writing  $u^k \partial_x u = \frac{1}{k+1} \partial_x (u^{k+1})$, we observe that this non-linear term essentially behaves as the derivative of a polynomial in $u$ which agrees with the classical assumption in the study of non-linear dispersive waves \cite{Abde,Saut}.\\

Finally, for $\eta>0$ fix, and moreover, for $n\in \mathbb{N}^{*}$ and $m=2,3$, the \emph{dissipative} effects are given by the  \emph{non-local} perturbation term  $\ds{\eta(\mathcal{H} \partial^{n}_x u +   \mathcal{H}_m u)}$. Here,  $\mathcal{H}$ denotes the Hilbert transform  defined  as:
\begin{equation}
\label{HIL}
\mathcal{H}(\varphi)(x)=p.v. \frac{1}{\pi} \int_{\mathbb{R}} \frac{\varphi(y)}{y-x} dy,
\end{equation}  where we have $\ds{\widehat{\mathcal{H}(\varphi)}(\xi)= i\,\sign(\xi) \widehat{\varphi}(\xi)}$. Moreover, the operator $\mathcal{H}_m $ is defined by the expression 
\begin{equation}
\mathcal{H}_m u = \left\{ \begin{array}{ll} \vspace{2mm}  
(-1)^{m-1} \partial^{m}_{x} u, & \,\, \text{if} \quad m=2,  \\
\mathcal{H} \partial^{m}_{x} u, &\,\, \text{if} \quad m=3. \end{array} \right. 
\end{equation}
When $n=1$, the term $\ds{\eta(\mathcal{H} \partial_x u +   \mathcal{H}_m u)}$ arises in physical models  and it describes the wave's instability in a stratified fluid (see the short explanation below equation (\ref{OST}) for more details). However, we will also consider  higher values of the parameter $n$ which, from the mathematical point of view, will play an interesting role in the spatial decaying properties of solutions.\\   

As already mentioned, equation (\ref{Equation}) is a generic model and its major interest bases on the fact that it  contains as a particular case  several \emph{relevant physical models}. In order to motivate the interest of equation (\ref{Equation}) let us examine the following examples. Let us set $n=1$: 
\begin{enumerate}
\item[$A)$] For $D= \partial^{2}_{x}$, where $p(\xi)=-\vert \xi\vert^2$, $k=1$ and $m=3$, the equation  (\ref{Equation}) deals with a non local perturbed version of the celebrated  \emph{Korteweg-de Vries} (KdV) equation \cite{KDV}. This equation, also known as the \emph{Ostrovsky, Stepanyams and Tsimring} (OST) equation: 
\begin{equation}\label{OST}
 \partial_t u + \partial^{3}_{x} u +u  \partial_x u + \eta \mathcal{H} (\partial_x u + \partial^{3}_{x} u) = 0,   
\end{equation}
describes the  radiational  instability of long non-linear waves in a stratified flow caused by internal wave radiation from a shear layer. The parameter $\eta>0$ represents the importance of amplification and damping relative to dispersion. The fourth term in equation represents amplification, while the fifth term in equation  denotes damping. For a more complete physical description we refer to  \cite{OST0,OST,OST1}.\\
\item[$B)$] For $D= \partial^{2}_{x}$, $k=2,3$ and $m=3$, the equation (\ref{Equation})  coincides with the generalized OST-equation:  
\begin{equation}\label{gOST}
 \partial_t u + \partial^{3}_{x} u +u^k  \partial_x u + \eta \mathcal{H} (\partial_x u + \partial^{3}_{x} u) = 0.   
\end{equation}
This model considers a stronger non-linear dynamics due to the term $u^k  \partial_x u$ with $k \geq 2$. The values $k=2$ and $k=3$ are relevant  from the physical point of view in the modelling of surface and volume water waves respectively \cite{Bona1}.  
\item[$C)$] For $D=\mathcal{H}\partial_x$, where $p(\xi)=\vert \xi \vert$, $k=1$ and $m=3$, the equation (\ref{Equation}) agrees with a non local perturbed version of the well-known  \emph{Benjamin-Ono} (BO) equation \cite{Benjamin}:  
\begin{equation}\label{BO}
\partial_t u + \mathcal{H} \partial^{2}_{x} u +u  \partial_x u + \eta \mathcal{H} (\partial_x u + \partial^{3}_{x} u) = 0.    
\end{equation}
This equation is a good approximate model for long-crested unidirectional waves at the interface of a two-layer system of incompressible inviscid fluids. Moreover, it gives an analogous model of the OST-equation (\ref{OST}) in deep stratified fluids \cite{Benjamin}.\\

\item[$D)$] For $D=\mathcal{H}\partial_x$, $k=1$ and $m=2$,  the equation (\ref{Equation}) becomes the \emph{Chen-Lee equation} (CL) which deals with the BO equation with another kind of non-local perturbation: 
 \begin{equation}\label{Chen-Lee}
 \partial_t u + \mathcal{H} \partial^{2}_{x} u +u  \partial_x u + \eta \mathcal{H} (\partial_x u )- \partial^{2}_{x} u = 0.    
 \end{equation}
Chen-Lee equation  was  introduced by H. H. Chen and Y. C. Lee in \cite{Chen1} to describe nonlinear dynamical models of plasma turbulence. See also \cite{Chen2} for more details.  
\item[$E)$] For  $D= (\mathcal{H} \partial_x)^{1+\alpha}$ with $0<\alpha <1$, where $p(\xi)=\vert \xi \vert^{1+\alpha}$, $k=1$ and $m=3$, the equation (\ref{Equation}) writes down as a non local perturbation (since we assume $\eta>0$) of the dispersive generalized BO equation:
\begin{equation}\label{dgBO}
\partial_t u + (\mathcal{H} \partial_x)^{1+\alpha} \partial_x  u +u  \partial_x u =0.    
\end{equation}
From the physical point of view, this  equations  models vorticity waves in the coastal zone \cite{SH}. On the other hand, from the mathematical point of view, this equation was studied in \bl{\cite{FonLinPon}}  as an  interesting \emph{intermediate} dispersive model between the BO equation (when $\alpha=0$) and the KdV equation (when $\alpha=1$). The parameter $0<\alpha<1$  measures the sharp  dispersive effects which are stronger than the one for the BO equation but  weaker than the one for the KdV equation.   
\end{enumerate}	
Concerning the mathematical study of these equations, as the local and global well-posedness (LWP and GWP respectively) and some previous results on  the spatial  decaying properties, there exists a large amount of works. So, we will give a short overview on the most recent results. \\

First, let us focus on the non-local perturbations of the KdV equation. For the \emph{OST-equation} (\ref{OST}), GWP was proved in $H^s(\R)$ with $s\geq 0$ in \cite{CarGaSa} and LWP was obtained in $H^{s}(\R)$ with $-3/2 \leq s <0$ in \cite{CuiZhao}. Moreover, the value $s=-3/2$ is the critical one for the LWP  in the Sobolev spaces. Thereafter, the average decay of solutions was derived in \cite{BorysAlvarez-2} using the weighted space $H^2\cap L^2((1+\vert \cdot \vert^2)dx)(\R)$. On the other hand, respect to the \emph{generalized OST-equation} (\ref{gOST}), only for the values $k=2$ and $k=3$, it was shown in \cite{CarScia} the LWP in $H^s(\R)$ for $s>0$ and the GWP in $L^2(\R)$. To the best of our knowledge, the well-posedness issues for $k\geq 4$ and the spatially decay properties for $k\geq 2$ have not been yet studied.\\

Now, let us concentrate on the non-local perturbations of the BO equation. The \emph{non-local perturbed BO equation} (\ref{BO}) was recently studied in \cite{FonPasRo} where the GWP was obtained in $H^s(\R)$  with $s>-3/2$. The value $s=-3/2$ seems to be critical for the well-posedness in the Sobolev spaces in the sens that the flow map data$-$solution for this  equation is not $\mathcal{C}^2$ from $H^s(\R)$ to $H^s(\R)$ for  $s <-3/2$.  Moreover, similar to the equation (\ref{OST}), the averaged decay of solutions was studied in the space $H^2\cap L^2((1+\vert \cdot \vert^2)dx)(\R)$. On the other hand, for the \emph{CL equation} (\ref{Chen-Lee}), the GWP was first  proved in \cite{PasRia2} for the periodic Sobolev spaces $H^s(\mathbb{T})$ with $s>-1/2$. Thereafter, this result was generalized to the non-periodic setting of the whole line $\R$ in \cite{PasRia}. Moreover, in this work it was also proved that the value $s=-1/2$ seems to the sharp  provided that  the flow map data$-$solution for this  equation is not $\mathcal{C}^3$ from $H^s(\R)$ to $H^s(\R)$ for   $s<-1/2$.  Finally, concerning the decay of solutions, always in \cite{PasRia} it was shown that solutions cannot have an averaged decay at infinity faster than  $1/ \vert x \vert^3$. More precisely, it is proven that if $u(t,x)$ is a solution of equation  (\ref{Chen-Lee}) which verifies $u \in \mathcal{C}([0,T], H^3\cap L^2(1+\vert \cdot \vert^6 dx)(\R))$ then we have $u(t,\cdot)=0$ for all $t\in [0,T]$.\\ 

Within the general framework of the equation (\ref{Equation}), the  aim of this paper is to give a better understanding of the spatial behavior of solutions for all these equations and other related models in the studying of stratified fluids. \\ 

Our methods are technically different with respect to the  works mentioned above. Indeed, these results are obtained through purely dispersive approaches based on  Strichartz-type estimates, smoothing effects and estimates in Bourgain-type spaces. Instead, following some ideas of a previous work \cite{CorJa}, using  the explicit definition in the frequency variable  and the inverse Fourier transform we derive some   sharp estimates (in the spatial variable) on the kernel associated to the linear part of (\ref{Equation}).\\

Kernels estimates seems to be a useful tool to study the spatial properties of solutions for equation (\ref{Equation}). First,  we derive  some  \emph{pointwise  decaying rates} of solutions which in certain cases are \emph{optimal}. Thereafter, combining the kernel estimates  with  some well-known tools of modern harmonic  analysis, as the Hardy–Littlewood maximal function operator and the  Muckenhoupt weights, we   study the \emph{average decay properties} of solutions in the improved setting of the weighted Lebesgue spaces $L^{p}_{w_\delta}(\R)$ for $1<p<+\infty$ (see formula (\ref{peso}) for a definition of these spaces).  Finally,  we are able to construct solutions of  with some \emph{pointwise growing rates}. To the best of our knowledge, these kind of solutions have not been considered in the previous works.\\

\textbf{Plan of the paper:}  in Section \ref{Sec:Kernel-Estimates} we derive some kernel estimates. Then, Section \ref{Sec:WP} is devoted to a first result on the well-posedness of equation (\ref{Equation}). In this framework, we study some pointwise and average decaying properties of solution in Sections \ref{Sec:PWD} and \ref{AVD} respectively. Finally, in Section \ref{PWG} we study some pointwise growing properties of solutions.  

\section{Statement of the results}
In all the results obtained in this paper,  we will observe that  the properties of solutions of equation (\ref{Equation}), as the well-posedness issues and the spatial behavior, deeply rely with parameters $m$ and $n$ given in the term $\ds{ \eta(\mathcal{H} \partial^{n}_x u +   \mathcal{H}_m u)}$. From now on we will assume that the parameter $n\in \mathbb{N}^{*}$ verifies 
\begin{equation}\label{Cond-n}
n \neq 5+ 4d, \quad \text{with}\quad d=0,1,2,3, \cdots.
\end{equation}
This condition on the parameter $n$ seems to be sharp to prove the well-posedness of equation (\ref{Equation}) in the classical Sobolev spaces $H^s(\R)$. Roughly speaking, solutions  can be written as  an explicit integral formulation (see formula (\ref{Equation-Int})) which involves a kernel depending on $m$ and $n$. For the values of $n$ which do not verify  (\ref{Cond-n}), the Fourier transform of the kernel behaves at high frequencies as an increasing exponential function (see the expression (\ref{FK1}) below) and then we loose any control on the well-known $H^s$-norm. Moreover, it is worth to emphasize that the condition (\ref{Cond-n}) is not too restrictive since all the physical models mentioned above are not concerned.     
\subsection{Global well-posedness}
To the best of our knowledge, the fairly  general equation (\ref{Equation}) has not been considered before in the literature and, in order to provide a more complete study of this equation, we give first the following result concerning some well-posedness issues in the classical framework of Sobolev spaces. 
\begin{TheoremeP}\label{Th-GWP} In equation (\ref{Equation}),  let the parameters $m=2,3$, $n \in \mathbb{N}^{*}$  which verifies  (\ref{Cond-n}) and  $k \in \mathbb{N}^{*}$. Moreover, for  $s>3/2$  let $u_{0} \in H^s(\R)$ be an initial datum. Then, the equation (\ref{Equation}) possess a unique classical solution 
 $$u \in \mathcal{C}([0,+\infty[, H^s(\R))\cap  \mathcal{C}^{1}(]0, +\infty[, \mathcal{C}^{\infty}(\R)).$$ 
\end{TheoremeP}
Some comments are in order. For the \emph{OST equation} (\ref{OST}) studied in \cite{CuiZhao}, the \emph{non-local perturbed BO equation} (\ref{BO}) treated in \cite{FonPasRo} and the \emph{CL equation} (\ref{Chen-Lee}) studied in \cite{PasRia},  this theorem recovers some well-known results obtained on the GWP and improves the regularity of solutions: here we have $\ds{u \in \mathcal{C}^{1}(]0, +\infty[, \mathcal{C}^{\infty}(\R))}$.  Moreover, for the \emph{generalized OST equation} (\ref{gOST}) studied in \cite{CarScia}, this theorem improves both  GWP and regularity of solutions for the higher non-linearities  $k\geq 4$. \\ 

It is worth mention that,  under some technical modifications, this theorem could be improved in two different ways. On one hand, for the particular case $k=1$, we are able to consider initial data belonging to the space $H^s(\R)$ with $s\geq 0$.  On the other hand, this result could be improved for a large class of initial data $u_0$ belonging to a Besov space $B^{s,q}_{2}(\R)$ with $2 \leq q \leq +\infty$. For a definition and more details on Besov spaces see the book \cite{BaChDan}. \\   

In the  framework of this theorem, we study now some decaying properties of solutions in the spatial variable.  
\subsection{Pointwise decaying properties}\label{SubSec:pointwise-decay}
Before to state our results we need to precise some notation. We denote $V \subset \R$ a neighborhood of the origin, and moreover, for $k \in \mathbb{N}$   we use the standard notation $C^{k}(V)$ for the functions $k-$ times differentiable in the neighborhood  $V$. \\

In our second result, we study the pointwise spatial decaying of solutions of the general equation (\ref{Equation}), provided that the initial datum $u_0$ verifies a pointwise decay. It is interesting to observe that  the parameter $n$ in the dissipative perturbation term: $\ds{\eta(\mathcal{H} \partial_x^{n} u +   \mathcal{H}_m u)}$, as well as the regularity properties at the origin  of the symbol $p(\xi)$ of the operator $D$ in the dispersive term : $\ds{D(\partial_x u)}$,  play a fundamental  role in the description of the spatial  decaying  of solutions.     
\begin{TheoremeP}\label{Th-Decay} Let $u_0 \in H^s(\R)$ (with $s>3/2$) be an initial datum and let $u(t,x)$ be the solution of equation (\ref{Equation}) given by Theorem \rd{\ref{Th-GWP}}.  Assume that for a parameter $\gamma>0$ the initial datum verifies 
\begin{equation}\label{Cond-datum}
 \vert u_0 (x)\vert \leq \frac{c}{\vert x \vert^{\gamma}}, \quad \vert x \vert \to +\infty.
\end{equation}
For $n \in \mathbb{N}^{*}$ given by (\ref{Cond-n}), if the symbol $p(\xi)$ verifies $p \in \mathcal{C}^{n-1}(V)$ then we have 
\begin{equation}\label{Decay1}
\vert u(t,x) \vert \leq   \frac{C_0}{\vert x \vert^{\min(\gamma,  n+1)}}, \quad \vert x \vert \to +\infty,   
\end{equation}
for a constant $C_0=C_0(u, t)>0$, depending on  the solution $u$ and the time $t>0$.
\end{TheoremeP}

Let us make the following remarks. In (\ref{Decay1}) we may observe that  the parameter $n$ \emph{controls} the decaying rate of solutions: the function $u(t,x)$ fulfills the decaying given by the initial datum  only if $\gamma \leq n+1$. But, if  the initial datum decays fast enough, when  $\gamma > n+1$, then the solution $u(t,x)$ does not mimic this decaying and it decays at infinity as $1 / \vert x \vert^{n+1}$. Moreover, we observe that the decaying properties given in (\ref{Decay1}) also depend of a \emph{equilibrium} between the dispersive and dissipative terms in equation (\ref{Equation}) in the sense that for higher values of parameter $n$ (in the dissipative term)  the symbol $p(\xi)$ (in the dispersive term) must be more regular at the origin. \\ 

From the physical point of view, the value $n=1$ is the most interesting  since the non-local perturbation  term $\ds{\eta(\mathcal{H} \partial_x u +   \mathcal{H}_m u) }$ (with $m=2$ or $m=3$) gives a good model of  long non-linear waves deformation in  stratified flows \cite{OST0, OST, OST1}. To illustrate the relevance of Theorem \ref{Th-Decay} in the studying of some  physical models,  let us mention the following examples. For the reader's convenience,  we will divide our study concerning   two relevant groups of equations.  
  
\begin{enumerate}
\item[$\bullet$]  \textbf{The \emph{KdV}-type models}. 

  For the \emph{OST-equation} (\ref{OST}), numerical studies  done in \cite{Feng} by B.F.  Feng  \&  T.  Kawahara shows that  for   every $ \eta>0$ there exists  a  family  of   solitary  waves  which \emph{experimentally} decay  as $\ds{1 / \vert x \vert^2}$  when $|x| \to + \infty$. Thus, setting $\gamma=2$ and  assuming that the initial datum verifies $\ds{\vert u_0 (x) \vert \leq c / \vert x \vert^2}$, then by  (\ref{Decay1}) we rigorously obtain the  decay rate of solutions
\begin{equation}\label{Dec2}
\vert u (t,x)\vert \lesssim \frac{1}{\vert x \vert^{2}}.    
\end{equation}
This decay rate was also exhibit in a previous work \cite{CorJa}. Moreover,  we observe that this decay rate is  also verified for the case of higher non-linearities in the setting \emph{generalized  OST-equation}   (\ref{gOST}).  \\

Recall that for $\eta=0$, the equation  (\ref{OST}) becomes the \emph{KdV} equation:

$$  \partial_t u + \partial^{3}_{x} u +u  \partial_x u=0.$$

For this equation,   T. Kato \cite{Kato} showed the  persistence problem in the setting of the space  $H^{2m}(\R) \cap L^{2}(|x|^{2m}, \,dx)$, with $m\in \mathbb{N}^{*}$.  Similarly, when $\eta=0$, the equations (\ref{gOST}) writes down as the \emph{generalizd KdV} equation:
$$ \partial_t u + \partial^{3}_{x} u +u^{k}  \partial_x u=0, \quad k \in \mathbb{N}^{*}.$$
In  \cite{NahasPonce},  J. Nahas \& G. Ponce obtained the same persistence results proved by T. Kato for this generalized version of the \emph{KdV} equation.  


\item[$\bullet$]  \textbf{The \emph{Benjamin-Ono}-type models.}

For the \emph{non-local perturbed BO} (\ref{BO}) and the \emph{CL equation} (\ref{Chen-Lee}), a second work \cite{Feng2} due to  B.F.  Feng  \&  T.  Kawahara shows \emph{numerically} that the solitary waves of these equations behave at infinity as $\ds{1 / \vert x \vert^2}$. Thus, always by  (\ref{Decay1}) we able to exhibit solutions of these equations which a explicit decay rate given in (\ref{Dec2}). Moreover, always in \cite{Feng2}, it is experimentally shown that if the perturbation parameter $\eta>0$ is large enough then the dispersive  term in equations (\ref{BO}) and (\ref{Chen-Lee}) can be  negligible (from the numerical point of view) and in this setting there exist solitary waves which behave at infinity as $\ds{1 / \vert x \vert}$. Thus, setting now $\gamma=1$, by (\ref{Decay1}) we get solutions of these equations with a pointwise decay of the form 
\begin{equation}\label{Dec1}
  \vert u (t,x) \vert \lesssim \frac{1}{\vert x \vert}.   
\end{equation} 
It is worth remark that  both equations (\ref{BO}) and (\ref{Chen-Lee}) agree with the \emph{Benjamin-Ono equation} when $\eta=0$:
$$ \partial_t u + \mathcal{H} \partial^{2}_{x} u +u  \partial_x u=0.$$
The spatial decaying properties for this equation are  also studied by  J. Nahas \& G. Ponce in \cite{NahasPonce}, where the following result is proven:  if $u_0 \in H^{2}(\R) \cap L^{2}(|x|^{4}, \,dx)$, then  the Cauchy problem for the \emph{Benjamin-Ono equation} equation is globally well-posed in the space $C([0,+ \infty[;  H^{2}(\R) \cap L^{2}(|x|^{4} \,dx))$.  \cite{NahasPonce}.  

\end{enumerate}

Although the main physical relevance of Theorem \ref{Th-Decay} is when $n=1$,  from the mathematical point of view  it is also  interesting to study the influence of high values  of the parameter $n$ in  the decaying behavior of solutions. As already mentioned, for the values  $n\geq 2$ the description of the  decay  of solutions becomes more \emph{complex} in the sense that it is also determined by the regularity properties of the symbol $p(\xi)$ at the origin $\xi=0$. Let us illustrate  this interesting phenom with some simple examples. For simplicity, we set $\gamma>0$ large enough, so we let the initial datum decay fast enough, and moreover we set $m=3$. 
\begin{enumerate}
\item[$\bullet$] For $n=2$, let us consider the following \emph{theoretical} non-local perturbation of the \emph{dispersive generalized BO equation} (\ref{dgBO}):  
\begin{equation*}
  \partial_t u + (\mathcal{H} \partial_x)^{1+\alpha} \partial_x  u  +u  \partial_x u + \eta(\mathcal{H} \partial^{2}_x u +   \mathcal{H} \partial^{3}_{x} u), \quad \text{with} \quad 0 < \alpha < 1.  \end{equation*}  
Here we have $p(\xi)=\vert \xi \vert^{1+\alpha}$ and then $p \in \mathcal(C)^{1}(V)$. Thus, by (\ref{Decay1}) the solutions have the following spatial behavior 
$$ \vert u(t,x) \vert \lesssim \frac{1}{\vert x \vert^3}, \quad  \vert x \vert \to +\infty.$$
\item[$\bullet$] For $n \geq 3$, let us consider the following \emph{theoretical} perturbed \emph{KdV} equation: 
\begin{equation*}
  \partial_t u + \partial^{3}_{x} u  +u  \partial_x u + \eta(\mathcal{H} \partial^{n}_x u + \mathcal{H} \partial^{3}_{x}u), \quad \text{with} \quad 0 < \alpha < 1.  \end{equation*}  
In this case we have $p(\xi)=-\vert \xi \vert^2$ and then $p \in \mathcal{C}^{\infty}(V)$. Thus, by (\ref{Decay1}) the solutions decay as follows: 
$$ \vert u(t,x) \vert \lesssim \frac{1}{\vert x \vert^{n+1}}, \quad  \vert x \vert \to +\infty.$$
\end{enumerate}

Now, it is natural to ask if the decay rates (\ref{Decay1})  are either optimal or they can be improved. In our third result, assuming  some technical restrictions on the parameters $m$ and $n$, we are able to answer these questions.  As was pointed out the in \cite{CorJa}, the zero-mean properties of the initial datum $u_0$ is the key tool to study these facts. \\

\begin{TheoremeP}\label{Th-Decay2}  Within the framework of Theorem \rd{\ref{Th-Decay}}, assume that $(m,n)\neq (2,1)$ and $(m,n)\neq (2, 2 d)$ with $d \in \mathbb{N}^{*}$. Moreover, let $0<\varepsilon \leq  1$.
\begin{enumerate}
\item[$1)$]  If the initial datum $u_0$ verifies:  \begin{equation*}
 \vert u_0 (x)\vert \leq \frac{c}{\vert x \vert^{n+1+\varepsilon}}, \quad \vert x \vert \to +\infty, \quad \text{and}\quad \int_{\R} u_0 (y) dy = 0,
\end{equation*} 
then the solution $u(t,x)$ of equation (\ref{Equation}) verifies:
\begin{equation}\label{Decay1Imp}
 \vert u (t,x)\vert \leq \frac{C_1}{\vert x \vert^{n+1+\varepsilon}},\quad \vert x \vert \to +\infty,    
\end{equation}
for a constant  $C_1=C_1(u, t)>0$. 
\item[$2)$] If the initial datum $u_0$ verifies: \begin{equation*}
 \vert u_0 (x)\vert \leq \frac{c}{\vert x \vert^{n+1+\varepsilon}}, \quad \vert x \vert \to +\infty, \quad \text{and}\quad \int_{\R} u_0 (y) dy \neq  0,
\end{equation*} 
then the solution $u(t,x)$ of equation (\ref{Equation}) verifies:
\begin{equation}\label{DecayOp}
C_2 \left\vert \int_{\R}u_0(y)dy \right\vert \frac{1}{\vert x \vert^{n+1}} \leq \vert u (t,x)\vert,\quad \vert x \vert \to +\infty,    
\end{equation}
for a constant  $C_2=C_2(t)>0$. 
\end{enumerate}
\end{TheoremeP}
In point $1)$, we may observe here that if the initial datum is a zero-mean function and if it decays fast enough (we have $\gamma=n+1+\varepsilon$) then the decay rate obtained in (\ref{Decay1}) is  improved in (\ref{Decay1Imp})  for $0<\varepsilon \leq 1$. To the best of our knowledge, the value $\varepsilon =1$ seems to be the maximal one to improve the  decay rates. This is due to the fact that the solutions of equation (\ref{Equation}) are written in a explicit integral formulation (\ref{Equation-Int}), where the spatial decay properties of the kernel eventually  block an improvement in the decaying of the solution for $\varepsilon >1$. On the other hand, in point $2)$, we remark that when the initial datum is not a zero-mean function then  the decay rate is  \emph{optimal}. Moreover, even if this datum is a fast-decaying function we have an instantaneous lost of persistence of solution $u(t,x)$ when $t>0$.\\

It is worth to emphasize  that  the additional conditions of the parameters $m$ and $n$ stated above are essentially technical and we refer the reader to Remark \ref{CondParam} in page \pageref{CondParam} for the details. However, this additional conditions are not too restrictive since most of the physical models mentioned in the introduction are considered in Theorem \ref{Th-Decay2}. Indeed, observe that  we can set the values $m=3$ and $n=1$   and then Theorem \ref{Th-Decay2} hols true for the
for the following  \emph{relevant physical models} that we also study considering the following groups of equations mentioned above. 

\begin{enumerate}
\item[$\bullet$] \textbf{The \emph{KdV}-type models}. 

For the equations  (\ref{OST}) and (\ref{gOST}),   we observe that if we consider a particular initial datum $u_0$ such that $\ds{u_{0}(x)= c / \vert x \vert^{\gamma}}$, for $\gamma>0$ large enough and for $\vert x \vert$ large enough;  then by point $2)$ we  obtain solutions  whit the \emph{sharp} asymptotic behavior: 
$$ \vert u(t,x) \vert \sim  \frac{1}{\vert x \vert^2},\quad \vert x \vert \to +\infty,$$ 
which agrees with the numerical results obtained in \cite{Feng} and \cite{Feng2}.  This optimal decaying results strongly differ to the classical \emph{KdV} equation, where  the solutions decay as fast as the initial datum \cite{Kato,NahasPonce}. This  different behavior  is caused by the effects of the perturbation term.

\item[$\bullet$] \textbf{The \emph{Benjamin-Ono}-type models}.	\\

For the classical \emph{Benjamin-Ono} equation, in \cite{NahasPonce} it is proven that the solutions of the \emph{Benjamin-Ono} equation  cannot decay at infinity faster than $1/ \vert x \vert^4$, whereas for the perturbed equation (\ref{BO}) a previous result, obtained by  G. Fonseca, R. Pastr\'an \& G. Rodr\'iguez-Blanco in \cite{FonPasRo}, proves  that the solutions of this equation cannot decay faster than $1/\vert x \vert^3$.  Here, we improve this last result result and we prove that, under the hypothesis of Theorem \ref{Th-Decay2}, the solutions have an optimal decay rate of the order $1 / \vert x \vert^2$. \\

Finally, let us mentions that  among these physical models  we only left open the case of the \emph{CL equation} (\ref{Chen-Lee}) which deals with the values $m=2$ and $n=1$ that are not included in this theorem.
\end{enumerate}	

\subsection{Average decaying properties} 
Our methods also allow us to study the average decay properties of solution $u(t,x)$. These  decay properties are characterized through the weighted Lebesgue  space which we introduce as follows: for the parameter $\gamma>0$ we introduce the weight
\begin{equation}\label{peso}
  w_\gamma (x)= \frac{1}{(1+\vert x \vert)^\gamma},   
\end{equation}  and for $1\leq p \leq +\infty$ we consider the weighted  Lebesgue  space
 $\ds{L^{p}_{w_\gamma}(\R)= L^p(w_\gamma \, dx)}$. The weighted  Lebesgue  spaces give us a fairly general framework to study different decaying properties of solutions of equation (\ref{Equation}). For the classical \emph{Lebesgue spaces} we have the (evident) embedding $L^p(\R) \subset L^{p}_{w_\gamma}(\R)$, but the weighted  Lebesgue  spaces also contains more \emph{sophisticate} functional spaces which characterize the average decaying of functions.  Denoting as $L^{p,\infty}(\R)$ a \emph{Lorentz space} (see the book \cite{DCh} for a complete study of these spaces)  we have the continuous embedding  $L^{p,\infty}(\R)\subset L^{p}_{w_\gamma}(\R) $. Moreover, for $1<r<p<+\infty$ we denote as $\dot{M}^{r,p}(\R)$ the \emph{homogeneous Morrey space} (see the Section $8$ of the book \cite{PGLR1} for a definition and some properties of these spaces). Then, for $0<1-r/p<\gamma$ we have the continuous embedding  $\dot{M}^{r,p}(\R)\subset L^{p}_{w_\gamma}(\R)$. In the setting of the weighted  Lebesgue  spaces we have the following result. 
\begin{TheoremeP}\label{Th-Averaged-Decay}  Let $u_0 \in H^s(\R)$ (with $s>3/2$) be an initial datum and let $u(t,x)$ be the solution of equation (\ref{Equation}) given by Theorem \rd{\ref{Th-GWP}}. For $1<p<+\infty$ and $0<\gamma <1$ assume that the initial datum verifies  $u_0 \in L^{p}_{w_\gamma}(\R)$. Then, for the parameter $\alpha >0$ given in (\ref{Alpha}) which only depends on $m=2,3$ and $n$ given by (\ref{Cond-n}), we have 
$$ u \in  L^{\infty}_{loc}\left(]0,+\infty[, L^{p}_{w_\gamma}(\R), t^{\alpha}\,dt\right). $$
\end{TheoremeP}  

Remark  that the have the continuous embedding $\ds{L^{2}(1+\vert x \vert^2)dx) \subset L^{2}_{w_\delta}(\R)}$. Thus, for $p=2$ this theorem improves some well-known results on the average decaying properties of solutions for \emph{non-local perturbed BO equation}  (\ref{BO}) studied  in  \cite{FonPasRo}, the \emph{CL-equation} (\ref{Chen-Lee}) treated in \cite{PasRia} and the  \emph{OST-equation}   (\ref{OST}) studied in  \cite{BorysAlvarez-2}. For this latter equation, due to the embedding $L^p(\R)\subset L^{p}_{w_\delta}(\R)$ for any $1<p<+\infty$, this theorem also improves a recent result on the average decaying properties given in \cite{CorJa}. Moreover, to our knowledge, this kind of results seems not be studied before for the \emph{generalized OST-equation}  (\ref{gOST}).   
\subsection{Pointwise growing properties} 
In all our previous results, we consider an initial datum $u_0$ with pointwise or average \emph{decaying} properties at infinity. However, it is also interesting to study the persistence problem of solutions for equation (\ref{Equation}) when the initial datum has some \emph{growing} properties at infinity. For the value $k=1$ in the non-linear term in (\ref{Equation}), we are able give a first result on the existence of solutions  which fulfill some pointwise growing properties (in the spatial variable) given by initial datum. \\

\begin{TheoremeP}\label{Th-Grow} Let $u_0 \in \dot{H}^{1}(\R)$ be an initial datum such that for $0<\gamma<1/2$ and for a constant $C_0>0$,  it verifies for all $x \in \R$: $$\vert u_0(x)\vert \leq C_0 (1+\vert x \vert)^{\gamma}.$$ Moreover, let $0<T<+\infty$.   There exists a constant $\delta=\delta(T)>0$ such that if $\ds{ \left\Vert  u_0 \right\Vert_{\dot{H}^{1}}+ C_0 < \delta}$,  then there exits a unique mild solution $u(t,x)$ of equation (\ref{Equation}) (with $k=1$) defined on the interval of time  $[0,T]$, such that  for all $x \in \R$ we have
$$ \vert u(t,x) \vert \leq  C (1+\vert x \vert)^{\gamma},$$ for a constant $C=C(u_0,u,t)>0$ depending on $u_0$, $u$ and $t$. 
\end{TheoremeP}

Let us make the following comments.  We observe first that this theorem does not come from the setting of Theorem \ref{Th-GWP} since, due to the well-known Sobolev embedding,  the assumption of the initial datum given in Theorem \ref{Th-GWP}: $u_0 \in H^s(\R)$ with $s>1/2$ implies  that $u_0 \in L^{\infty}(\R)$ which is not coherent with the growing properties assumed above. In this theorem, we assume instead $u_0 \in \dot{H}^{1}(\R)$ and this hypothesis is essentially technical. However, it is worth to remark that this hypothesis is coherent with the growing properties assumed. A simple example of an initial datum verifying all the  hypothesis in Theorem \ref{Th-Grow} is given by 
\begin{equation*}
u_0(x)= \left\{ \begin{array}{cc}\vspace{2mm}
C_0(1+x)^{\gamma},& x >0, \\
0, & x \leq 0. \end{array}\right.       
\end{equation*}
Here, as  $0<\gamma<1/2$ it is easy to see that we have $\frac{d}{dx}u_0 \in L^2(\R)$.\\

To close this section, let us mention that from now on in the following computations, the generic constants $C_\eta>0$ and $c_\eta>0$ may change in each line but they  only depend on the fixed parameter  $\eta>0$ given in the dissipative perturbation term in equation (\ref{Equation}).

\section{Kernel estimates}\label{Sec:Kernel-Estimates}
Remark first that the equation (\ref{Equation}) can be written as the equivalent integral formulation 
\begin{equation}\label{Equation-Int}
u(t,x)= K_{m,n}(t,\cdot)\ast u_0(x)+ \int_{0}^{t} K_{m,n}(t-\tau) \ast (u^{k}\partial_x u)(\tau, x) d\tau, 
\end{equation} where, for $t>0$, $m=2,3$, $n\in \mathbb{N}^{*}$, and moreover, for the symbol $p(\xi)$ given in formula (\ref{def-P}), the kernel $K_{m,n}(t,x)$ is defined in the Fourier variable as follows:
\begin{equation}\label{Def-K-Fou}
\widehat{K_{m,n}}(t,\xi)=e^{-i p(\xi) \xi t -\eta (i^{n+1} \vert \xi \vert \xi^{n-1} +\vert \xi \vert^m)t}. 
\end{equation} To make the notation more convenient let us introduce the function 
\begin{equation}\label{def-f2n}
\varphi_{m,n}(\xi)= -\eta (i^{n+1} \vert \xi \vert \xi^{n-1} + \vert \xi \vert^m)= \left\{\begin{array}{ll}\vspace{2mm} 
\eta (i^{n+1} \xi^n -(-\xi)^m), \quad  \xi <0, \\
-\eta(i^{n+1} \xi^n + \xi^m), \quad  \xi \geq 0.\\
\end{array}\right.
\end{equation} With this notation write $\ds{\widehat{K_{m,n}}(t,\xi)=e^{-i p(\xi) \xi t + \varphi_{m,n}(\xi)t}}$, hence, as the symbol $p(\xi)$ is a real-valued function we obtain 
$$ \vert \widehat{K_{m,n}}(t,\xi)\vert = \vert e^{\varphi_{m,n}(\xi)\, t} \vert.$$
In this expression we are interesting in the behavior of the quantity $\ds{\varphi_{m,n}(\xi)}$ which comes from the dissipative perturbation term in equation (\ref{Equation}). A simple calculation shows that  for $m=2,3$ and for $n=5+4d$, with $d \in \mathbb{N}$, we have $\ds{\varphi_{m,n}(\xi) = \eta(\vert \xi \vert^n - \vert \xi \vert^m) t}$, then  we get \begin{equation}\label{FK1}
\vert e^{\varphi_{m,n}(\xi)\, t} \vert  \sim e^{\eta \vert \xi \vert^n \,t}, \quad \vert \xi \vert \to  +\infty,   
\end{equation}
and thus,  for those values of $n$ we loose any control on the function $\ds{\widehat{K_{m,n}}(t,\xi)}$. On the other hand, observe that for the values of $n$ which verify the condition (\ref{Cond-n}), \emph{i.e.}, $n \neq 5+4d$, for  $n$ even we have the identity  $\ds{\varphi_{m,n}(\xi) = \eta(i \vert \xi \vert \xi^{n-1}-\vert \xi \vert^m)}$, hence we obtain 
\begin{equation}\label{FK2}
 \vert e^{\varphi_{m,n}(\xi)\, t} \vert  \sim e^{-\eta \vert \xi \vert^m\,t}, \quad \vert \xi \vert \to  +\infty,    
\end{equation}
moreover, for $n$ odd we have the identity $\ds{\varphi_{m,n}(\xi) = -\eta(\vert \xi \vert^n + \vert \xi \vert^m)}$ and then we get 
\begin{equation}\label{FK3}
 \vert e^{\varphi_{m,n}(\xi)\, t} \vert  \sim e^{-\eta(\vert \xi \vert^n + \vert \xi \vert^m)\,t}, \quad \vert \xi \vert \to  +\infty.    
\end{equation}
In conclusion, when $n$ verifies (\ref{Cond-n})  the function $\ds{\widehat{K_{m,n}}(t,\xi)}$ has good decaying properties and the key idea is to use them to obtain sharp estimates on the kernel in the spatial variable. Notice that by  (\ref{Cond-n})  the kernel $K_{m,n}(t,x)$ writes down as the convergent integral: 
\begin{equation}\label{def-K}
K_{m,n}(t,x)=\int_{\R} e^{2 \pi i x \xi}  \, e^{-i p(\xi) \xi t +\varphi_{m,n}(\xi)t} \, d\xi.
\end{equation}
In the following result  we study some spatial decaying properties of the kernel  $\ds{K_{m,n}(t,x)}$, which will be fundamental  in the next sections. As mentioned in Section \ref{SubSec:pointwise-decay},  the regularity properties at the origin of the symbol $p(\xi)$   plays an important role  in this study  and they vary as long as the parameter $n$ take different values. Thus, for the sake of clarity, we will consider first the values $n=1,2$ and then  the values $n \geq 3$.
\begin{Proposition}\label{Prop-Kernel} Let  $m=2,3$ and $n \in \mathbb{N}^{*}$ which verifies (\ref{Cond-n}).  There exist two constants $C_\eta>0$ and $c_\eta>0$, depending only on $\eta>0$,   such that for all $t>0$ and for all $x \in \R$ it verifies:
	\begin{enumerate}
		\item[$1)$] For $n=1, 2$, if the symbol $p(\xi)$ verifies $p \in \mathcal{C}^{n-1}(V)$ then   we have:   $\ds{\vert K_{m,n}(t,x) \vert \leq C_\eta \, \frac{e^{c_\eta  t }}{\eta ^{1/m} t^{1/m}}  \frac{1}{1 + \vert x \vert^{n+1}}}$. 
		\item[$2)$] For  $n \geq 3$,  if   symbol $p(\xi)$  verifies $p \in \mathcal{C}^{n-1}(V)$ then we have: \begin{enumerate}
			\item[$2.1)$] For $n$ even	($n=2d$, with $d \in \mathbb{N}$ and $d\geq 2$): $\ds{\vert K_{m,n}(t,x) \vert \leq c_\eta \, \frac{e^{3 \eta t}}{\eta^{1/m} t^{1/m}} \frac{1}{1+ \vert x \vert^{n+1-\ell}}}$. 
			\item[$2.2)$] For $n$ odd	($n=3+4d$, with $d \in \mathbb{N}$): $ \ds{\vert K_{m,n}(t,x) \vert \leq c_\eta\, \frac{e^{3 \eta t }}{\eta^{1/n} t^{1/n}} \frac{1}{1+ \vert x \vert^{n+1-\ell}}}$. 
		\end{enumerate}	 
	\end{enumerate}	
\end{Proposition}

\pv  We start  writing  $$\ds{ K_{m,n}(t,x)  =\int_{\xi<0} e^{2 \pi i x \xi} e^{-i p(\xi) \xi t   \  +\varphi_{m,n}(\xi) t}\, d\xi + \int_{\xi>0} e^{2 \pi i x \xi} e^{-i p(\xi) \xi t   \  +\varphi_{m,n}(\xi) t}\, d\xi}.$$ The key idea to study the pointwise spatial behavior of the kernel $K_{m,n}(t,x)$ essentially bases on the procedure: first, in each term of the right-hand side in  this identity, for $x\neq 0$ we multiply and we divide by $2\pi \, i x$ to obtain:
\begin{equation}\label{Prod1}
\begin{split}
K_{m,n}(t,x)= & \frac{1}{2\pi \, i x} \int_{\xi<0} (2\pi \, i x)e^{2 \pi i x \xi} e^{-i p(\xi) \xi t   \  +\varphi_{m,n}(\xi) t}\, d\xi + \frac{1}{2\pi \, i x} \int_{\xi>0}(2\pi \, i x) e^{2 \pi i x \xi} e^{-i p(\xi) \xi t   \  +\varphi_{m,n}(\xi) t}\, d\xi \\
=& \frac{1}{2 \pi i x } \int_{\xi<0} \partial_\xi (e^{2 \pi i x \xi})  e^{-i p(\xi) \xi t   \  +\varphi_{m,n}(\xi) t}\, d\xi + \frac{1}{2 \pi i x } \int_{\xi>0} \partial_\xi (e^{2 \pi i x \xi})  e^{-i p(\xi) \xi t   \  +\varphi_{m,n}(\xi) t}\, d\xi.
\end{split}
\end{equation} Thereafter,  integrating by parts each term we write  
\begin{equation}\label{Prod2}
\begin{split}
K_{m,n}(t,x)= &  \frac{1}{2 \pi i x } 
	\left( e^{2 \pi i x \xi} e^{-i p(\xi) \xi t   \  +\varphi_{2,n}(\xi) t}  \Big|^{0}_{-\infty}  - \int_{\xi <0} e^{2 \pi i x \xi} \partial_{\xi} (e^{-i p(\xi) \xi t   \  +\varphi_{m,n}(\xi) t}) \, d \xi  \right)\\
	&  +\frac{1}{2 \pi i x } \left(e^{2 \pi i x \xi} e^{-i p(\xi) \xi t   \  +\varphi_{m,n}(\xi) t}
	\Big|^{+\infty}_{0} -  \int_{\xi >0} e^{2 \pi i x \xi} \partial_{\xi} (e^{-i p(\xi) \xi t   \  +\varphi_{m,n}(\xi) t}) \, d \xi  \right) \\
	=& \frac{1}{2 \pi i x } \left( 1 - \int_{\xi <0} e^{2 \pi i x \xi} \partial_{\xi} (e^{-i p(\xi) \xi t   \  +\varphi_{m,n}(\xi) t}) \, d \xi  \right) \\
	& +\frac{1}{2 \pi i x } \left(  -1  -  \int_{\xi >0} e^{2 \pi i x \xi} \partial_{\xi} (e^{-i p(\xi) \xi t   \  +\varphi_{m,n}(\xi) t}) \, d \xi  \right) \\
	=& - \frac{1}{2 \pi i x} \left( \int_{\xi <0} e^{2 \pi i x \xi} \partial_{\xi} (e^{-i p(\xi) \xi t   \  +\varphi_{m,n}(\xi) t}) \, d \xi +\int_{\xi >0} e^{2 \pi i x \xi} \partial_{\xi} (e^{-i p(\xi) \xi t   \  +\varphi_{m,n}(\xi) t}) \, d \xi  \right).
	\end{split} 
\end{equation}
Repeating the computations done in (\ref{Prod1}) and (\ref{Prod2}) we get 
\begin{equation*}
\begin{split}
 	K_{m,n}(t,x)= &- \frac{1}{(2 \pi i x)^2} \left( e^{2 \pi i x \xi} \partial_{\xi} (e^{-i p(\xi) \xi t   \  +\varphi_{m,n}(\xi) t})  \Big|^{0}_{-\infty} -  \int_{\xi <0} e^{2 \pi i x \xi} \,  \partial^{2}_{\xi} (e^{- i p(\xi) \xi t + \varphi_{m,n}(\xi) t}) \, d \xi \right)\\
	&  - \frac{1}{(2 \pi i x)^2} \left(  e^{2 \pi i x \xi} \partial_{\xi} (e^{-i p(\xi) \xi t   \  +\varphi_{m,n}(\xi) t})    \Big|^{+\infty}_{0} -  \int_{\xi >0} e^{2 \pi i x \xi} \,  \partial^{2}_{\xi} (e^{- i p(\xi) \xi t + \varphi_{m,n}(\xi) t}) \, d \xi \right).   
\end{split}  
\end{equation*}
Now, remark that we have 
\begin{equation}\label{DER1}
e^{2 \pi i x \xi} \partial_{\xi} (e^{-i p(\xi) \xi t   \  +\varphi_{m,n}(\xi) t}) = e^{2 \pi i x \xi}   e^{- i p(\xi) \xi t + \varphi_{m,n}(\xi) t} (-i p'(\xi) \xi t -i p(\xi) t +\varphi_{m,n}^{'}(\xi)t),    
\end{equation}
hence, by the good decaying properties of the function $\ds{e^{\varphi_{m,n}(\xi) t}}$ when $\vert \xi \vert \to +\infty$ (see the formulas (\ref{FK2}) and (\ref{FK3})) and moreover,   as by (\ref{def-f2n}) we have $\varphi_{m,n}(0)=0$, then we get 
$$ e^{2 \pi i x \xi}   e^{- i p(\xi) \xi t + \varphi_{m,n}(\xi) t} (-i p'(\xi) \xi -i p(\xi)  +\varphi_{m,n}^{'}(\xi))t  \,  \Big|^{0}_{-\infty} = -i p(0^{-})t +\varphi^{'}_{m,n}(0^{-})t,$$
and
$$e^{2 \pi i x \xi}   e^{- i p(\xi) \xi t + \varphi_{m,n}(\xi) t} (-i p'(\xi) \xi -i p(\xi) +\varphi_{m,n}^{'}(\xi))t  \,  \Big|^{+\infty}_{0} = i p(0^{+})t- \varphi^{'}_{m,n}(0^{+})t,$$ where, for a function $f(\xi)$ we use the standard notation  $\ds{\lim_{\xi \to 0^{-}}f (\xi)=f(0^{-})}$ and  $\ds{\lim_{\xi \to 0^{+}}f (\xi)=f(0^{+})}$. \\

Thus we can write 
\begin{equation}\label{eqAUX1}
\begin{split}
K_{m,n}(t,x)=& - \frac{1}{(2 \pi i x)^2}\left( ip(0^{-})t-ip(0^{+})t -\varphi^{'}_{m,n}(0^{-}) t+ \varphi^{'}_{m,n}(0^{+})t\right)  \\
 &+\frac{1}{(2 \pi i x)^2}  \left( \int_{\xi <0} e^{2 \pi i x \xi} \,  \partial^{2}_{\xi} (e^{- i p(\xi) \xi t + \varphi_{m,n}(\xi) t}) \, d \xi+ \int_{\xi >0} e^{2 \pi i x \xi} \,  \partial^{2}_{\xi} (e^{- i p(\xi) \xi t + \varphi_{m,n}(\xi) t}) \, d \xi \right).  
\end{split}
\end{equation}
At this point, in order to make a clearer exposition, we will distinguish the following cases of the parameter $n$.\\

\textbf{The case $n=1$.} Recall that in this case  we  assume $p\in \mathcal{C}(V)$ and then we get $p(0^{-})-p(0^{+})=0$. Moreover, the function  $\varphi^{'}_{m,n}(\xi)$ writes down as: 
\begin{equation}\label{eq02}
\varphi^{'}_{m,n}(\xi)= \left\{    \begin{array}{ll} \vspace{2mm} 
\eta(i^{n+1} n \xi^{n-1} + m(-\xi)^{m-1}), & \, \xi <0, \\ 
-\eta(i^{n+1} n \xi^{n-1} + m\xi^{m-1}), & \, \xi >0, 
\end{array} \right.
\end{equation}
hence, for $n=1$ and $m=2,3$  we have $-\varphi_{m,1}^{'}(0^{-})+\varphi_{m,1}^{'}(0^{+})=2 \eta$. Then, getting back to (\ref{eqAUX1}) we can write:
\begin{eqnarray}\label{aux1} \nonumber
K_{m,1}(t,x) &=&- \frac{2 \eta t }{(2 \pi i x)^2} + \frac{1}{(2 \pi i x)^2}  \left(  \int_{\xi <0} e^{2 \pi i x \xi} \,  \partial^{2}_{\xi} (e^{- i p(\xi) \xi t + \varphi_{m,1}(\xi) t}) \, d \xi + \int_{\xi >0} e^{2 \pi i x \xi} \,  \partial^{2}_{\xi} (e^{- i p(\xi) \xi t + \varphi_{m,1}(\xi) t}) \, d \xi  \right)\\
&=& - \frac{2 \eta t }{(2 \pi i x)^2} + I_1.
\end{eqnarray}
Once the term $\ds{-\frac{2 \eta t }{(2 \pi i x)^2}}$ appears, this fact suggests the estimate $\ds{\vert K_{m,1}(t,x)\vert \lesssim 1 / \vert x \vert^2}$. Indeed, we study the term $I_1$ above, where, applying the computations done in (\ref{Prod1}) and (\ref{Prod2}) we get 
\begin{equation*}
\begin{split}
I_1=& \frac{1}{(2\pi i x)^3} \left(e^{2 \pi i x \xi} \,  \partial^{2}_{\xi} (e^{- i p(\xi) \xi t + \varphi_{m,1}(\xi) t})  \Big|^{0}_{-\infty}  - \int_{\xi <0} e^{2 \pi i x \xi} \partial^{3}_{\xi} (e^{-i p(\xi) \xi t   \  +\varphi_{m,1}(\xi) t}) \, d \xi  \right)\\
&  +\frac{1}{(2 \pi i x)^3} \left(e^{2 \pi i x \xi} \,  \partial^{2}_{\xi} (e^{- i p(\xi) \xi t + \varphi_{m,1}(\xi) t})
	\Big|^{+\infty}_{0} -  \int_{\xi >0} e^{2 \pi i x \xi} \partial^{3}_{\xi} (e^{-i p(\xi) \xi t   \  +\varphi_{m,1}(\xi) t}) \, d \xi  \right) \\
=& \frac{1}{(2\pi i x)^3} (I_{1,1}+I_{1,2}). 	
\end{split}
\end{equation*}
Now, by the good decaying properties of the function $e^{\varphi_{m,1}(\xi) t}$ when $\vert \xi \vert \to +\infty$ (see the expression (\ref{FK3}) with $n=1$) and following the same computations done in   Lemma $3.1$ of \cite{CorJa}  we have $\vert I_{1,1}+ I_{1,2} \vert \leq C_\eta \, e^{c_\eta t}$. Then we obtain the following estimate \begin{equation}\label{I1}
\vert I_1 \vert \leq C_\eta \frac{e^{c_\eta t}}{\vert x \vert^3}.    
\end{equation}
Thus,  by (\ref{aux1}) and (\ref{I1}), for $\vert x \vert$ large enough we can write 
\begin{equation}\label{aux2}
\vert K_{m,1}(t,x) \vert \leq c \frac{2 \eta t}{\vert x \vert^2}+ \vert I_1 \vert  \leq c \frac{2 \eta t}{\vert x \vert^2}+ C_\eta \frac{e^{c_\eta t}}{\vert x \vert^3}\leq C_\eta \frac{e^{c_\eta t}}{\vert x \vert^2}+ C_\eta \frac{e^{2 \eta t}}{\vert x \vert^3} \leq C_\eta \frac{e^{ c_\eta t}}{\vert x \vert^2}.      
\end{equation}
Until now  we have estimated the kernel $K_{m,1}(t,x)$ for $\vert x \vert$ sufficiently  large and it remains to obtain an estimate also valid for $\vert x \vert$ small. For this   we have write 
\begin{equation}\label{eqAUX4}
  \vert K_{m,1}(t,x) \vert \leq  \Vert K_{m,1}(t,\cdot) \Vert_{L^{\infty}} \leq \Vert \widehat{K_{m,1}} (t,\cdot) \Vert_{L^1},  
\end{equation} where, by (\ref{def-f2n}) the quantity $\ds{\Vert \widehat{K_{m,1}} (t,\cdot) \Vert_{L^1}}$ is estimated as follows: \begin{eqnarray*}
\Vert \widehat{K_{m,1}} (t,\cdot) \Vert_{L^1} &\leq &   \int_{\R} \left\vert e^{\varphi_{m,1}(\xi) t} \right\vert  d \xi \leq \int_{\R} e^{\eta(\vert \xi \vert -\vert\xi\vert^m)t} d \xi \leq  \int_{\vert \xi \vert \leq 2^{2^{1/(m-1)}}} e^{\eta(\vert \xi \vert -\vert\xi\vert^m)t} d \xi \\
& & + \int_{\vert \xi \vert \geq  2^{1/(m-1)}}  e^{\eta(\vert \xi \vert -\vert\xi\vert^m)t} d \xi \\
&\leq &  \int_{\vert \xi \vert \leq 2^{1/(m-1)}} e^{\eta \vert \xi \vert t} d \xi + \int_{\vert \xi \vert \geq 2^{1/(m+1)}} e^{-\eta \frac{\vert \xi\vert^m}{2}t}  d \xi \\
 &\leq & c e^{2^{1/(m-1)} \eta t} + c \frac{1}{(\eta t )^{1/m}} \leq C_\eta \frac{e^{c_\eta t}}{(\eta t)^{1/m}}.
\end{eqnarray*}
Hence, for all $x\in \R$ we get $\ds{\vert K_{m,1}(t,x)\vert \leq  C_\eta \frac{e^{c_\eta t}}{(\eta t)^{1/m}}}$. Finally, gathering this estimate and the estimate given in (\ref{aux2}) we obtain 
\begin{equation}\label{Res1}
\vert K_{m,1}(t,x)\vert \leq C_\eta \frac{e^{c_\eta t}}{\eta^{1/m}\, t^{1/m}}\frac{1}{1+\vert x \vert^2}.     
\end{equation}
\textbf{The case $n=2$.} Recall that in this case we assume $p\in \mathcal{C}^{1}(V)$, hence, in particular we have  $p(0^{-})-p(0^{+})=0$. Moreover,  by (\ref{eq02}) for $n=2$ and $m=2,3$ we have  $\varphi_{m,2}^{'} (0^{-})=\varphi_{m,2}^{'}(0^{+})=0$. Then,   getting back to (\ref{eqAUX1}) we get 
\begin{equation*}
 K_{m,2}(t,x)= \frac{1}{(2 \pi i x)^2}  \left(  \int_{\xi <0} e^{2 \pi i x \xi} \,  \partial^{2}_{\xi} (e^{- i p(\xi) \xi t + \varphi_{m,2}(\xi) t}) \, d \xi + \int_{\xi >0} e^{2 \pi i x \xi} \,  \partial^{2}_{\xi} (e^{- i p(\xi) \xi t + \varphi_{m,2}(\xi) t}) \, d \xi  \right).    
\end{equation*}
Here we observe that  we can continue with the same process and we apply  computations done in (\ref{Prod1}) and (\ref{Prod2}) to obtain 
\begin{equation*}
\begin{split}
K_{m,2}(t,x)
=& \frac{1}{(2 \pi i x)^3}  \left(  e^{2 \pi i x \xi} \partial^{2}_{\xi} (e^{- i p(\xi) \xi t + \varphi_{m,2}(\xi) t}) \Big|^{0}_{-\infty} -  \int_{\xi <0} e^{2 \pi i x \xi} \,  \partial^{3}_{\xi} (e^{- i p(\xi) \xi t + \varphi_{m,2}(\xi) t}) \, d \xi \right)\\
&  +\frac{1}{(2 \pi i x)^3}  \left(  e^{2 \pi i x \xi} \partial^{2}_{\xi} (e^{- i p(\xi) \xi t + \varphi_{m,2}(\xi) t})  \Big|^{+\infty}_{0}-  \int_{\xi >0} e^{2 \pi i x \xi} \,  \partial^{3}_{\xi} (e^{- i p(\xi) \xi t + \varphi_{m,2}(\xi) t}) \, d \xi \right).
\end{split}
\end{equation*}
In order to study these expressions, remark first that we have
\begin{equation}\label{der2}
\begin{split}
 \partial^{2}_{\xi} (e^{- i p(\xi) \xi t + \varphi_{m,2}(\xi) t}) =& e^{- i p(\xi) \xi t + \varphi_{m,2}(\xi) t}  \left[\left(-i p'(\xi) \xi t -i p(\xi)t+ \varphi^{'}_{m,2}(\xi)t \right)^2  \right. \\
 & \left. -i p''(\xi) \xi t-2i p'(\xi)t + \varphi^{''}_{m,2}(\xi)t \right],    
\end{split}
\end{equation}
where the function $\ds{\varphi^{''}_{m,n}(\xi)}$ writes down as: 
\begin{equation}\label{2der}
\varphi^{''}_{m,n}(\xi)= \left\{    \begin{array}{ll} \vspace{2mm} 
\eta \left(i^{n+1}n(n-1)\xi^{n-2}-m(m-1)(-\xi)^{m-2}\right), & \, \xi <0, \\ 
-\eta \left(i^{n+1}n(n-1)\xi^{n-2}+m(m-1)\xi^{m-2}\right), & \, \xi >0. \\
\end{array} \right.      
\end{equation}
Recalling that by (\ref{eq02}) we have  $\varphi_{m,2}^{'}(0^{-})=\varphi_{m,2}^{'}(0^{+})=0$ and moreover, remarking that by  (\ref{2der}) we have $\varphi_{m,2}^{''}(0^{-})= - \eta(2i+c_m)$ and $\varphi_{m,2}^{''}(0^{+})= - \eta(-2i+c_m)$, where $c_m=2$ if $m=2$, and $c_m=0$ if $m=3$,  then, by the good decaying properties of the function $e^{\varphi_{m,2}(\xi) t}$ when $\vert \xi \vert \to +\infty$ (see the expression (\ref{FK2}) with $n=2$)   we get:  
$$  e^{2 \pi i x \xi} \partial^{2}_{\xi} (e^{- i p(\xi) \xi t + \varphi_{m,2} (\xi) t}) \Big|^{0}_{-\infty} =   -2i p'(0^{-})t - \eta(2i+c_m)t,$$ and 
$$ e^{2 \pi i x \xi} \partial^{2}_{\xi} (e^{- i p(\xi) \xi t + \varphi_{m,2} (\xi) t})  \Big|^{+\infty}_{0}= +2 i p'(0^{+})t +  \eta (-2i +c_m)t,$$ and then we can write
\begin{eqnarray*}
	K_{m,2}(t,x)&=&  \frac{1}{(2 \pi i x)^3} ( - 2 i  p'(0^{-})t + 2 i p'(0^{+}) t  - 4 i \eta t  )  \\
	& & - \frac{1}{(2 \pi i x)^3}  \left(  \int_{\xi <0} e^{2 \pi i x \xi} \,  \partial^{3}_{\xi} (e^{- i p(\xi) \xi t + \varphi_{m,2}(\xi) t}) \, d \xi +  \int_{\xi >0} e^{2 \pi i x \xi} \,  \partial^{3}_{\xi} (e^{- i p(\xi) \xi t + \varphi_{m,2}(\xi) t}) \, d \xi \right).
\end{eqnarray*}
But, recalling that we have $p\in \mathcal{C}^{1}(V)$ we finally obtain  
\begin{eqnarray}\label{Estim2} \nonumber
	K_{m,2}(t,x)&=&  \frac{- 4 i \eta t }{(2 \pi i x)^3}  - \frac{1}{(2 \pi i x)^3}  \left(  \int_{\xi <0} e^{2 \pi i x \xi} \,  \partial^{3}_{\xi} (e^{- i p(\xi) \xi t + \varphi_{m,2}(\xi) t}) \, d \xi +  \int_{\xi >0} e^{2 \pi i x \xi} \,  \partial^{3}_{\xi} (e^{- i p(\xi) \xi t + \varphi_{m,2}(\xi) t}) \, d \xi \right)\\ 
	&=& \frac{- 4 i \eta t }{(2 \pi i x)^3} + I_2.
\end{eqnarray}
We must study now the term $I_2$. By the good decaying properties of the function $e^{\varphi_{m,2}(\xi)t}$ (see always the expression (\ref{FK2}) with $n=2$) and moreover,   following the same computations performed  for the term $I_1$ in (\ref{I1})   the term $I_2$ is estimated as follows: 
\begin{equation}\label{I2}
\vert I_2 \vert \leq C_\eta \frac{e^{c_\eta t}}{\vert x \vert^4}.     
\end{equation} Then, for $\vert x \vert$ large enough  we can write 
\begin{equation}\label{eqAUX3}
\vert K_{m,2}(t,x)\vert \leq c \frac{\eta t }{\vert x \vert^3}+ C_\eta \frac{e^{c_\eta t}}{\vert x \vert^4} \leq C_\eta \frac{e^{c_\eta t }}{\vert x \vert^3}+ C_\eta \frac{e^{c_\eta t }}{\vert x \vert^4}\leq C_\eta \frac{e^{c_\eta t }}{\vert x \vert^3}. 
\end{equation}
On the other hand, by estimate (\ref{eqAUX4}) we have $\ds{\vert K_{m,2}(t,x)\vert \leq \widehat{K_{m,2}} (t,\cdot)\Vert_{L^1}}$, where, by (\ref{def-f2n}) we write 
\begin{eqnarray*}
& &\Vert \widehat{K_{m,2}} (t,\cdot) \Vert_{L^1} =  \int_{\R} \left\vert  e^{- i p(\xi) \xi t + \varphi_{m,2}(\xi) t}  \right\vert  d \xi = \int_{\R} \left\vert  e^{- i p(\xi) \xi t -\eta(-i \vert \xi \vert \xi + \vert \xi \vert^m \vert)  t}  \right\vert    d \xi \leq \int_{\R} e^{-\eta \vert\xi\vert^m \, t} d \xi\\
& \leq  &  \frac{c}{(\eta t )^{1/m}} \leq C_\eta \frac{e^{c_\eta t}}{(\eta t)^{1/m}},
\end{eqnarray*} and then, for all $x \in \R$ we get $\ds{\vert K_{m,2}(t,x)\vert \leq   C_\eta \frac{e^{c_\eta t}}{(\eta t)^{1/m}}}$.  Finally, gathering this estimate and estimate (\ref{eqAUX3}) we have 
\begin{equation}\label{Res2}
\vert K_{m,2}(t,x)\vert \leq C_\eta \frac{e^{c_\eta t}}{\eta^{1/m}\, t^{1/m}}\frac{1}{1+\vert x \vert^3}.  
\end{equation}
At this point, by estimates (\ref{Res1}) and (\ref{Res2}) we have proven the point $1)$ in Proposition \ref{Prop-Kernel}. \\  

\textbf{The case $n\geq 3$.} The computations follow the same ideas performed in the previous cases ($n=1$ and $n=2$). Recall that in this case  we assume $p \in \mathcal{C}^{n-1}(V)$.  In particular we have $p \in \mathcal{C}(V)$ hence we get $p(0^{-})-p(0^{+})=0$. Moreover, by (\ref{eq02}) with $n\geq 3$ and $m=2,3$, we have $\varphi^{'}_{m,n}(0^{-})=\varphi^{'}_{m,n}(0^{+})=0$. Then, getting back to (\ref{eqAUX1}) we obtain 
$$ K_{m,n}(t,x)= \frac{1}{(2 \pi i x)^2}  \left( \int_{\xi <0} e^{2 \pi i x \xi} \,  \partial^{2}_{\xi} (e^{- i p(\xi) \xi t + \varphi_{m,n}(\xi) t}) \, d \xi+ \int_{\xi >0} e^{2 \pi i x \xi} \,  \partial^{2}_{\xi} (e^{- i p(\xi) \xi t + \varphi_{m,n}(\xi) t}) \, d \xi \right).$$
Thereafter, applying the computations done in (\ref{Prod1}) and (\ref{Prod2}) we get 
\begin{equation*}
\begin{split}
K_{m,n}(t,x) =& \frac{1}{(2 \pi i x)^3}  \left(  e^{2 \pi i x \xi} \partial^{2}_{\xi} (e^{- i p(\xi) \xi t + \varphi_{m,n}(\xi) t}) \Big|^{0}_{-\infty} -  \int_{\xi <0} e^{2 \pi i x \xi} \,  \partial^{3}_{\xi} (e^{- i p(\xi) \xi t + \varphi_{m,n}(\xi) t}) \, d \xi \right)\\
&  +\frac{1}{(2 \pi i x)^3}  \left(  e^{2 \pi i x \xi} \partial^{2}_{\xi} (e^{- i p(\xi) \xi t + \varphi_{m,n}(\xi) t})  \Big|^{+\infty}_{0}-  \int_{\xi >0} e^{2 \pi i x \xi} \,  \partial^{3}_{\xi} (e^{- i p(\xi) \xi t + \varphi_{m,n}(\xi) t}) \, d \xi \right).
\end{split}
\end{equation*} In this expression, by identity (\ref{der2}), the fact that $e^{\varphi_{m,n}(\xi)t}$ is a fast decaying function when $\vert \xi \vert \to +\infty$ (see the expression (\ref{FK2}) for $n$ even and the expression (\ref{FK3}) for $n$ odd) and moreover, since by (\ref{2der}) we have $\varphi^{''}_{m,n}(0^{-})= -\eta c_m$ and $\varphi^{''}_{m,n}(0^{+})=- \eta c_m$, with $c_m=2$ if $m=2$, and $c_m=0$ if $m=3$, we get  
$$  e^{2 \pi i x \xi} \partial^{2}_{\xi} (e^{- i p(\xi) \xi t + \varphi_{m,2} (\xi) t}) \Big|^{0}_{-\infty} =   -2i p'(0^{-})t - \eta c_m\,t,$$ and 
$$ e^{2 \pi i x \xi} \partial^{2}_{\xi} (e^{- i p(\xi) \xi t + \varphi_{m,2} (\xi) t})  \Big|^{+\infty}_{0}= +2 i p'(0^{+})t +  \eta c_m\,t.$$
Thus, we can write 
\begin{equation*}
\begin{split}
 K_{m,n}(t,x) =&  \frac{1}{(2 \pi \, i x)^3}  \left( 2 i t ( -p'(0^{-}) + p'(0^{+})) \right)\\
 &+  \frac{1}{(2 \pi i x)^3}  \left(   \int_{\xi <0} e^{2 \pi i x \xi} \,  \partial^{3}_{\xi} (e^{- i p(\xi) \xi t + \varphi_{m,n}(\xi) t}) \, d \xi  +  \int_{\xi >0} e^{2 \pi i x \xi} \,  \partial^{3}_{\xi} (e^{- i p(\xi) \xi t + \varphi_{m,n}(\xi) t}) \, d \xi  \right).
\end{split}
\end{equation*}
At this point recall that we have $p\in \mathcal{C}^{1}(V)$ (since $p\in \mathcal{C}^{n-1}(V)$ and $n\geq 3$) and then $-p'(0^{-}) + p'(0^{+})=0$. Then we obtain
\begin{equation}\label{eqAUX6}
 K_{m,n}(t,x) = \frac{1}{(2 \pi i x)^3}  \left(   \int_{\xi <0} e^{2 \pi i x \xi} \,  \partial^{3}_{\xi} (e^{- i p(\xi) \xi t + \varphi_{m,n}(\xi) t}) \, d \xi  +  \int_{\xi >0} e^{2 \pi i x \xi} \,  \partial^{3}_{\xi} (e^{- i p(\xi) \xi t + \varphi_{m,n}(\xi) t}) \, d \xi  \right).  
\end{equation}
At this point we observe that   we can apply the computations done in (\ref{Prod1}) and (\ref{Prod2}) 
iteratively  until to obtain the identity
\begin{equation}\label{eqAUX9}
\begin{split}
K_{m,n}(t,x) =& \frac{1}{(2 \pi i x)^{n+1}}  \left(  e^{2 \pi i x \xi} \partial^{n}_{\xi} (e^{- i p(\xi) \xi t + \varphi_{m,n}(\xi) t}) \Big|^{0}_{-\infty} -  \int_{\xi <0} e^{2 \pi i x \xi} \,  \partial^{n+1}_{\xi} (e^{- i p(\xi) \xi t + \varphi_{m,n}(\xi) t}) \, d \xi \right)\\
&  +\frac{1}{(2 \pi i x)^{n+1}}  \left(  e^{2 \pi i x \xi} \partial^{n}_{\xi} (e^{- i p(\xi) \xi t + \varphi_{m,n}(\xi) t})  \Big|^{+\infty}_{0}-  \int_{\xi >0} e^{2 \pi i x \xi} \,  \partial^{n+1}_{\xi} (e^{- i p(\xi) \xi t + \varphi_{m,n}(\xi) t}) \, d \xi \right). 
\end{split}    
\end{equation}
Here, as we may observe  in identities (\ref{DER1}) and (\ref{der2}), the expression $\ds{\partial^{n}_{\xi} (e^{- i p(\xi) \xi t + \varphi_{m,n}(\xi) t})}$ computes down as 
\begin{equation}\label{dern}
\partial^{n}_{\xi} (e^{- i p(\xi) \xi t + \varphi_{m,n}(\xi) t}) = e^{- i p(\xi) \xi t + \varphi_{m,n}(\xi) t} \left[ g_n(\xi, p(\xi), \varphi_{m,n}(\xi),t) -n\, i p^{(n-1)}(\xi)t + \varphi^{(n)}_{m,n}(\xi)t\right],    
\end{equation}
where $\ds{g_n(\xi, p(\xi),\varphi_{m,n}(\xi), t)}$ is a polynomial of degree $n$ which depends  on $\xi$, the derivatives $p^{(k)}(\xi)$ and $\varphi^{(k)}_{m,n}(\xi)$ for $k=0, \cdots, n-2$, and $t$, and moreover, it verifies $\ds{g_n(0, p(0),\varphi_{m,n}(0), t)=0}$.  On the other hand, the function $\ds{\varphi
^{(n)}_{m,n}(\xi)}$ computes down as: 
\begin{equation}\label{nder}
\varphi^{(n)}_{m,n}(\xi)= \left\{    \begin{array}{ll} \vspace{2mm} 
\eta(i^{n+1}n !+c_m), & \, \xi <0, \\ 
-\eta(i^{n+1}n ! + c_m), & \, \xi >0, \\
\end{array} \right.      
\end{equation} here, when  $n=3$ we have  $c_m=0$ if $m=2$ and $c_m=6$ if $m=3$, and moreover, when $n>3$ we have $c_m=0$ for $m=2,3$. \\ 

Thus, by (\ref{dern}) and (\ref{nder}) we obtain 
$$  e^{2 \pi i x \xi} \partial^{n}_{\xi} (e^{- i p(\xi) \xi t + \varphi_{m,n} (\xi) t}) \Big|^{0}_{-\infty} =   -n\,i p^{(n-1)}(0^{-})t  + \eta (i^{n+1}n !+c_m)\,t,$$ and 
$$ e^{2 \pi i x \xi} \partial^{n}_{\xi} (e^{- i p(\xi) \xi t + \varphi_{m,2} (\xi) t})  \Big|^{+\infty}_{0}= n\, i p^{n-1}(0^{+})t +  \eta (i^{n+1}n ! + c_m)\,t,$$
and getting back to (\ref{eqAUX9}), as we have $p\in \mathcal{C}^{n-1}(V)$ we are able to write 
\begin{eqnarray}\label{Estim3} \nonumber
& & K_{m,n}(t,x)= \frac{1}{(2 \pi i x)^{n+1}} \left(-n\,i p^{(n-1)}(0^{-})t+ n\, i p^{n-1}(0^{+})t + 2 \eta (i^{n+1} n!+c_m) t \right) \\ \nonumber
& &- \frac{1}{(2 \pi i x)^{n+1}} \left(  \int_{\xi <0} e^{2 \pi i x \xi} \,  \partial^{n+1}_{\xi} (e^{- i p(\xi) \xi t + \varphi_{m,n}(\xi) t}) \, d \xi +  \int_{\xi >0} e^{2 \pi i x \xi} \,  \partial^{n+1}_{\xi} (e^{- i p(\xi) \xi t + \varphi_{m,n}(\xi) t}) \, d \xi \right)\\ \nonumber
& =& \frac{c_\eta\, t}{(2 \pi i x)^{n+1}} -  \frac{1}{(2 \pi i x)^{n+1}} \left(  \int_{\xi <0} e^{2 \pi i x \xi} \,  \partial^{n+1}_{\xi} (e^{- i p(\xi) \xi t + \varphi_{m,n}(\xi) t}) \, d \xi +  \int_{\xi >0} e^{2 \pi i x \xi} \,  \partial^{n+1}_{\xi} (e^{- i p(\xi) \xi t + \varphi_{m,n}(\xi) t}) \, d \xi \right)\\
& =& \frac{c_\eta\, t}{(2 \pi i x)^{n+1}} + I_n. 
\end{eqnarray}
Thereafter, always by the good decaying properties of the function $e^{\varphi_{m,n}(\xi) t}$ (see the expression (\ref{FK2}) for $n$ even and the expression (\ref{FK3}) for $n$ odd) and moreover, following the computations done for the term $I_1$ in (\ref{I1})  we have the estimate 
\begin{equation}\label{I3}
\vert I_n \vert \leq C_\eta \frac{e^{c_\eta \, t}}{\vert x \vert^{n+2}}.   
\end{equation} 
Then, for $\vert x \vert$ large enough we obtain
\begin{equation}\label{eqAUX10}
\vert K_{m,n}(t,x) \vert \leq  \frac{c_\eta \, t}{\vert x \vert^{n+1}}+ C_\eta \frac{e^{c_\eta \, t}}{\vert x \vert^{n+2}} \leq C_{\eta} \frac{e^{c_\eta \, t}}{\vert x \vert^{n+1}} + C_\eta \frac{e^{c_\eta \, t}}{\vert x \vert^{n+2}} \leq C_\eta \frac{e^{c_\eta \, t}}{\vert x \vert^{n+1}}.    
\end{equation} 
On the other hand, by estimate (\ref{aux2}) we have $\ds{\vert K_{m,n}(t,x)\vert \leq \Vert \widehat{K_{m,n}}(t,\xi)\Vert_{L^1}}$, where the quantity $\ds{\Vert \widehat{K_{m,n}}(t,\xi)\Vert_{L^1}}$ is estimated as follows. For $n$ even, by (\ref{FK2}) we have 
$$ \Vert \widehat{K_{m,n}} (t,\cdot) \Vert_{L^1}   \leq \int_{\R} e^{-\eta \vert\xi\vert^m \, t} d \xi  \leq    \frac{c}{(\eta t )^{1/m}} \leq C_\eta \frac{e^{c_\eta t}}{(\eta t)^{1/m}}.$$ 
In this case,  for all $x$ we have $\ds{\vert K_{m,n}(t,x)\vert \leq C_\eta \frac{e^{c_\eta t}}{(\eta t)^{1/m}}}$ and with this estimate and estimate (\ref{eqAUX10}) we obtain $\ds{\vert K_{m,n}(t,x)\vert \leq  C_\eta \frac{e^{c_\eta t}}{(\eta^{1/m}\, t^{1/m}} \frac{1}{1+\vert x \vert^{n+1}}}$,  which proves the point $2.1)$ in Proposition \ref{Prop-Kernel}. On the other hand, for $n$ odd, by (\ref{FK3}) we have  
$$ \Vert \widehat{K_{m,n}} (t,\cdot) \Vert_{L^1}   \leq \int_{\R} e^{-\eta \vert\xi\vert^n \, t} d \xi  \leq    \frac{c}{(\eta t )^{1/n}} \leq C_\eta \frac{e^{c_\eta t}}{(\eta t)^{1/n}}.$$
Here  for  all $x$ we have $\ds{\vert K_{m,n}(t,x)\vert \leq C_\eta \frac{e^{c_\eta t}}{(\eta t)^{1/n}}}$ and then, by this estimate and by estimate estimate (\ref{eqAUX10}) we get $\ds{\vert K_{m,n}(t,x)\vert \leq  C_\eta \frac{e^{c_\eta t}}{\eta^{1/n}\, t^{1/n}} \frac{1}{1+\vert x \vert^{n+1}}}$,  which proves the point $2.2)$ in Proposition \ref{Prop-Kernel}. This proposition in now proven.  \finpv

Finally, in order to simplify the notation, for  $m=2,3$ and $n \in \mathbb{N}$ such that  (\ref{Cond-n}) is verified, let us define the parameter $\alpha>0$:
\begin{equation}\label{Alpha} 
\alpha= \left\{  \begin{array}{ll}\vspace{2mm} 
1/m, & \text{if} \, \, n=1 \,\, \text{or $n$ even:} \,\, n=2d \,\, \text{with} \, \, d \in \mathbb{N}\, \, \text{and}\, \, d\geq 2, \\
1/n, & \text{if $n$ odd:}\,\, n=3+4d, \,\, \text{with}\, \, d \in \mathbb{N}, 
\end{array} 
\right.
\end{equation}
hence,  we have $0<\alpha \leq 1/2$.   With this parameter, and the  estimates of the kernel $K_{m,n}(t,x)$ given in Proposition \rd{\ref{Prop-Kernel}}, from now on  we write the following unified kernel estimate:
\begin{equation}\label{Kernel-estimates}
\vert K_{m,n}(t,x) \vert \leq C_\eta\, \frac{e^{c_\eta t }}{ t^{\alpha}} \frac{1}{1+ \vert x \vert^{n+1}}.   
\end{equation}

\section{Global well-posedness: proof of Theorem \rd{\ref{Th-GWP}}}\label{Sec:WP}
 We will start by the following local well-posedness result.
\begin{Proposition}\label{LocalExistence} Let $s>3/2$ and let $u_0\in H^s(\R)$ be an initial datum. There exists a time $0<T_0 <+\infty$ and a function $u\in \mathcal{C}([0, T_0], H^s(\R))$, which is the  unique solution of equation (\ref{Equation-Int}).
\end{Proposition}
\pv For a time $0<T<+\infty$ (which we will set small enough) we consider the Banach space $\mathcal{C}([0, T], H^s(\R))$ with the usual norm $\ds{\Vert u \Vert_{T}=\sup_{0\leq t \leq T} \Vert u(t,\cdot)\Vert_{H^s}}$. \\

For the first term in the right-hand side of (\ref{Equation-Int}) we have $K_{m,n}(t,\cdot)\ast u_0 \in \mathcal{C}([0, T], H^s(\R))$. Indeed, remark first that by (\ref{FK2}) and $(\ref{FK3})$ there exists a constant $c=c(m,n,\eta)>0$, such that for all $t>0$ and all $\xi \in \R$ we have  $\ds{\vert \widehat{K_{m,n}}(t,\xi) \vert \leq c}$. Then, as $u_0 \in H^s(\R)$  we get 
\begin{equation}\label{estim-lin1}
 \sup_{0\leq t\leq T} \Vert K_{m,n}(t,\cdot) \ast u_0 \Vert_{H^s} \leq c \Vert u_0 \Vert_{H^s}.   
\end{equation}
It remains to prove the continuity of the quantity $\Vert K_{m,n}(t,\cdot) \ast u_0 \Vert_{H^s}$ on $[0,T]$. By convergence dominated we get directly $\ds{\lim_{t \to 0^{+}} \Vert K_{m,n}(t,\cdot) \ast u_0 -u_0 \Vert_{H^s}=0}$. Moreover, we have the following technical lemma:
\begin{Lemme}\label{Lemma-tech-2} Let  $\varepsilon >0$  and let $s_1 >0 $. Then, there exists a constant $c_1>0$, which depends on $s_1$, $\varepsilon$ and the parameters $m,n$, such that  for all $\varepsilon< t_1, t_2$, and for all $\psi  \in H^s(\R)$ we have:
	$$ \Vert K_{m,n}(t_1, \cdot)\ast \psi - K_{m,n}(t_2, \cdot)\ast \psi \Vert_{H^{s+s_1}} \leq c_1 \vert t_1 - t_2 \vert^{1/2} \Vert \psi \Vert_{H^s}.$$
\end{Lemme}	   
\pv Recall that by (\ref{Def-K-Fou}) and (\ref{def-f2n}) we have $\ds{\widehat{K_{m,n}}(t,\xi)=e^{-i p(\xi)\xi t + \varphi_{m,n}(\xi) t }}$. Then we can write 
\begin{equation*}
\begin{split}
 &\Vert K_{m,n}(t_1, \cdot)\ast \varphi - K_{m,n}(t_2, \cdot)\ast \varphi \Vert^{2}_{H^{s+s_1}}  = \int_{\R} (1+\vert \xi \vert^2)^{s+s_1} \vert \widehat{K_{m,n}}(t_1, \xi)- \widehat{K_{m,n}}(t_2, \xi)\vert^2 \vert \widehat{\psi}(\xi)\vert^2 d\, \xi  \\
 =& \int_{\R} (1+\vert \xi \vert^2)^{s+s_1} \vert e^{(-i p(\xi) \xi+\varphi_{m,n}(\xi) )t_2} \vert^{2} \vert e^{(-i p(\xi) \xi  +\varphi_{m,n}(\xi))(t_1-t_2)} -1\vert^2  \vert \widehat{\psi}(\xi)\vert^2 d\, \xi= (a),
\end{split}    
\end{equation*}
 where we must study the quantity $\ds{\vert e^{(-i p(\xi) \xi  +\varphi_{m,n}(\xi))(t_1-t_2)} -1\vert^2}$. We write 
 $$ \vert e^{(-i p(\xi) \xi  +\varphi_{m,n}(\xi))(t_1-t_2)} -1\vert^2= \vert e^{(-i p(\xi) \xi  +\varphi_{m,n}(\xi))(t_1-t_2)} -1\vert \, \vert e^{(-i p(\xi) \xi  +\varphi_{m,n}(\xi))(t_1-t_2)} -1\vert.$$
 Recall that by (\ref{FK2}) and (\ref{FK3}) the quantity $\ds{\vert e^{(-i p(\xi) \xi  +\varphi_{m,n}(\xi))(t_1-t_2)} -1\vert}$ is uniformly bounded and then we have 
$$ \vert e^{(-i p(\xi) \xi  +\varphi_{m,n}(\xi))(t_1-t_2)} -1\vert^2 \leq c \vert e^{(-i p(\xi) \xi  +\varphi_{m,n}(\xi))(t_1-t_2)} -1\vert.$$
Now, by  the mean value theorem in the temporal variable, the definition of $\varphi_{m,n}(\xi)$ given in (\ref{def-f2n}), and moreover,  by the estimate (\ref{estim-symbol}) on the symbol $p(\xi)$ we obtain:
$$ \vert e^{(-i p(\xi) \xi  +\varphi_{m,n}(\xi))(t_1-t_2)} -1\vert  \leq  c \vert -i p(\xi) \xi  + \varphi_{m,n}(\xi) \vert \vert t_1-t_2 \vert \leq  c (\vert \xi \vert^{\sigma+1} + \eta \vert \xi \vert^n + \eta \vert \xi \vert^{m} ) \vert t_1 - t_2 \vert.$$
Then we have $$ \vert e^{(-i p(\xi) \xi  +\varphi_{m,n}(\xi))(t_1-t_2)} -1\vert^2 \leq  c (\vert \xi \vert^{\sigma+1} + \eta \vert \xi \vert^n + \eta \vert \xi \vert^{m} ) \vert t_1 - t_2 \vert. $$
With this estimate, we get back to identity $(a)$ to write 
\begin{eqnarray*}
(a)&\leq & c \vert t_1 - t_2 \vert \int_{\R} (1+\vert \xi \vert^2)^{s+s_1} \vert e^{(-i p(\xi) \xi  +\varphi_{m,n}(\xi))t_2} \vert^{2} (\vert \xi \vert^{\sigma+1} + \eta \vert \xi \vert^n + \eta \vert \xi \vert^{m} ) \vert \widehat{\psi}(\xi)\vert^2 d\, \xi\\
&\leq & c \vert t_1 - t_2 \vert \int_{\R} (1+\vert \xi \vert^2)^{s+s_1} \vert e^{\varphi_{m,n}(\xi) \, t_2} \vert^{2} (\vert \xi \vert^{\sigma+1} + \eta \vert \xi \vert^n + \eta \vert \xi \vert^{m} ) \vert \widehat{\psi}(\xi)\vert^2 d\, \xi\\
&\leq & c \vert t_1 - t_2 \vert \, \underbrace{\sup_{\xi \in \R} \left( (1+\vert \xi \vert^2)^{s_1} (\vert \xi \vert^{\sigma+1} + \eta \vert \xi \vert^n + \eta \vert \xi \vert^{m} ) \vert e^{\varphi_{m,n}(\xi) t_2}\vert^{2} \right)}_{(b)} \Vert \psi \Vert^{2}_{H^s}.
\end{eqnarray*} At this point, recall that by (\ref{FK2}) and (\ref{FK3}), and moreover, as have  $t_2 > \varepsilon$ then  we obtain   
$$ (b) \leq  c \,  \sup_{\xi \in \R} \left( (1+\vert \xi \vert^2)^{s_1} (\vert \xi \vert^{\sigma+1} + \eta \vert \xi \vert^n + \eta \vert \xi \vert^{m} ) \max( e^{-2 \eta (\vert \xi \vert^n + \vert \xi \vert^m) \varepsilon}, e^{-2 \eta\vert \xi \vert^m \varepsilon} \right)= c_1 <+\infty, $$ hence the desired estimate follows.  \finpv
In this lemma we set $s_1=0$ and $\psi =u_0$, hence we obtain $K_{m,n}(t,\cdot)\ast u_0 \in \mathcal{C}(]0, T], H^s(\R))$. Thus, we finally have $K_{m,n}(t,\cdot)\ast u_0 \in \mathcal{C}[0, T], H^s(\R))$.  \\

We study now the second term in the right-hand side of (\ref{Equation-Int}). For this term have the following estimate: 
\begin{equation}\label{estim-nonlin1}
\sup_{0\leq t\leq T} \left \Vert  \int_{0}^{t} K_{m,n}(t-\tau, \cdot)\ast (u^k \partial_x u)(\tau, \cdot) d \tau \right\Vert_{H^s} \leq  C_{\eta} \, e^{c_\eta \, T}  T^{1-\alpha} \,  \Vert u \Vert^{k+1}_{T},  
\end{equation} where the parameter $0<\alpha \leq 1/2$ is given in (\ref{Alpha}). Indeed, for $0<t\leq T$ fix we write 
\begin{eqnarray*}
& & 	\left \Vert  \int_{0}^{t} K_{m,n}(t-\tau, \cdot)\ast (u^k \partial_x u)(\tau, \cdot) d \tau \right\Vert_{H^s} \leq \int_{0}^{t} \left\Vert K_{m,n}(t-\tau,\cdot)\ast (u^k\, \partial_x u)(\tau, \cdot) \right\Vert_{H^s} d\tau \\
&\leq & \frac{1}{k+1} \int_{0}^{t} \left\Vert K_{m,n}(t-\tau, \cdot) \ast \partial_{x} (u^{k+1})(\tau, \cdot) \right\Vert_{H^s} d \tau \leq  \frac{c}{k+1} \int_{0}^{t} \left\Vert K_{m,n}(t-\tau, \cdot) \ast u^{k+1}(\tau, \cdot) \right\Vert_{H^{s+1}} d \tau. 
\end{eqnarray*}	
At this point,  we need the following technical lemma:
\begin{Lemme}\label{Lemma-tech-1} Let $0<\alpha\leq 1/2$ be the parameter given in (\ref{Alpha}) and let $s_1, s_2 \in \R$. The following estimates follows: 
\begin{enumerate}
\item[$1)$] For   all $\psi \in H^{s_1}(\R)$ we have $\ds{\Vert K_{m,n}(t,\cdot) \ast \psi \Vert_{H^{s_1+s_2}} \leq C_{\eta, s_2} \frac{e^{c_{\eta\, s_2}\, t}}{t^{\alpha \,s_2}} \Vert \psi \Vert_{H^{s_1}}}$.
\item[$2)$] Moreover, all $\psi \in \dot{H}^{s_1}(\R)$ we have $\ds{\Vert K_{m,n}(t,\cdot) \ast \psi \Vert_{\dot{H}^{s_1+s_2}} \leq C^{'}_{\eta, s_2} \frac{e^{c^{'}_{\eta,s_2}\,   t}}{t^{\alpha \,s_2}} \Vert \psi \Vert_{\dot{H}^{s_1}}}$. 
\end{enumerate}
\end{Lemme}	
\pv The proof of points $1)$ and $2)$ essentially follows the same lines so it is sufficient to detail the computations for the point $1)$. We write 
\begin{eqnarray*}
& & \Vert K_{m,n}(t,\cdot)\ast \psi \Vert^{2}_{H^{s_1+s_2}} = \int_{\Rn} (1+\vert \xi \vert^2)^{s_1+s_2} \vert \widehat{K_{m,n}}(t,\xi) \vert^2 \vert \widehat{\psi} (\xi) \vert^2 \, d \xi\\
&= & \int_{\R} (1+\vert \xi \vert^2)^{s_2} \vert \widehat{K_{m,n}}(t,\xi) \vert^2 (1+\vert \xi \vert^2)^{s_1} \vert \widehat{\psi} (\xi) \vert^2 \, d \xi \leq \underbrace{ \left( \sup_{\xi \in \Rn}  (\vert 1+\vert \xi \vert^2)^{s_2} \vert \widehat{K_{m,n}}(t,\xi) \vert^2 \right)}_{(b)}\, \Vert \psi \Vert^{2}_{H^{s_1}}, 
\end{eqnarray*} 
where we must estimate the quantity $(b)$. For this we will consider the following cases of the parameters $m$ and $n$. For $n$ even, by  (\ref{FK2})  for all $t>0$ and for all $\xi \in \R$ we can write 
\begin{eqnarray*}
& & (\vert 1+\vert \xi \vert^2)^{s_2} \vert \widehat{K_{m,n}}(t,\xi) \vert^2  \leq  C_{s_2} \left( \vert \widehat{K_{m,n}}(t,\xi) \vert^2 + \vert \xi \vert^{2\, s_2} \vert \widehat{K_{m,n}}(t,\xi) \vert^2 \right)\leq   C_{\eta, s_2}\left( 1 + \vert \xi \vert^{2\, s_2} e^{-2 \eta \vert \xi \vert^{m} t} \right)\\
&\leq & C_{\eta, s_2} \left( 1+  \frac{\vert (\eta t)^{1/m} \xi \vert^{2 s_2}}{(\eta t)^{2 s_2 / m}} e^{-2  \vert (\eta t)^{1/m} \xi \vert^{m}} \right) \leq C_{\eta, s_2} \left( 1+ \frac{1}{(\eta t)^{2 s_2 / m}} \right)\leq C_{\eta, s_2} \frac{1+ (\eta t)^{2 s_2 / m}}{(\eta t)^{2 s_2 / m}} \leq C_{\eta, s_2}\frac{e^{c_{\eta, s_2}\, t }}{t^{2 s_2 / m}}. 
\end{eqnarray*} 
Then, for $n$ even  we have $\ds{(b) \leq C_{\eta, s_2}\frac{e^{c_{\eta, s_2}\, t }}{t^{2 s_2 / m}}}$, and thus we get 
\begin{equation}\label{eqAUX12}
\Vert K_{m,n}(t,\cdot)\ast \psi \Vert_{H^{s_1+s_2}} \leq C_{\eta, s_2}\frac{e^{c_{\eta, s_2}\, t }}{t^{ s_2 / m}}  \Vert \psi \Vert_{H^{s_1}}.   
\end{equation}
Now, for $n$ odd, by  (\ref{FK3}) and following the same estimates above, if $n \leq m$  we have $\ds{(b) \leq C_{\eta, s_2}\frac{e^{c_{\eta, s_2}\, t }}{t^{2 s_2 / n}}}$, and then we obtain: 
\begin{equation}\label{eqAUX13}
\Vert K_{m,n}(t,\cdot)\ast \psi \Vert_{H^{s_1+s_2}} \leq C_{\eta, s_2}\frac{e^{c_{\eta, s_2}\, t }}{t^{ s_2 / n}}  \Vert \psi \Vert_{H^{s_1}}.    
\end{equation} Moreover,  if $n>m$ then we have $\ds{(b) \leq C_{\eta, s_2}\frac{e^{c_{\eta, s_2}\, t }}{t^{2 s_2 / n}}}$, and thus we get the same estimate (\ref{eqAUX12}). Thus, by estimates (\ref{eqAUX12}) and (\ref{eqAUX13}), and moreover, by definition of the parameter $\alpha$ in (\ref{Alpha}) we obtain the estimate stated in point $1)$.  \finpv
In the setting of this lemma, we set the parameters $s_1=s$, $s_2=1$ and the function $\psi= u^{k+1}$. Moreover, as $s>1/2$ (since we have $s>3/2$) by the Sobolev  product laws  we can write 
\begin{eqnarray*}
& & \frac{c}{k+1} \int_{0}^{t} \left\Vert K_{m,n}(t-\tau, \cdot) \ast u^{k+1}(\tau, \cdot) \right\Vert_{H^{s+1}} d \tau \leq \frac{C_\eta}{k+1} \int_{0}^{t} \frac{e^{ c_\eta (t-\tau)}}{  (t-\tau)^{ \alpha}} \Vert u^{k+1}(\tau, \cdot) \Vert_{H^{s}} d \tau \\
&\leq & C_\eta \frac{e^{c_\eta\, T}}{(k+1)} \int_{0}^{t} \frac{\Vert u^{k+1}(\tau, \cdot) \Vert_{H^s}}{(t-\tau)^{\alpha}} d \tau \leq C_\eta \frac{e^{c_\eta\, T}}{(k+1)} \int_{0}^{t} \frac{\Vert u(\tau, \cdot) \Vert^{k+1}_{H^s}}{(t-\tau)^{\alpha}} d \tau \\
& \leq &C_\eta \frac{e^{c_\eta\, T}}{(k+1)}  \left(\sup_{0\leq t \leq T} \Vert u(\tau, \cdot) \Vert_{H^s} \right)^{k+1}\, \left(\int_{0}^{t} \frac{d \tau }{(t-\tau)^{\alpha}}\right)  \leq C_\eta \frac{e^{c_\eta\, T}}{(k+1)}  \Vert u \Vert^{k+1}_{T}\, t^{1-\alpha}\leq C_\eta e^{c_\eta\, T}  \Vert u \Vert^{k+1}_{T}\, t^{1-\alpha}, 
\end{eqnarray*} 
hence we obtain (\ref{estim-nonlin1}). Once we dispose of estimates (\ref{estim-lin1}) and (\ref{estim-nonlin1}), we set a time $0<T_0 <+\infty$ small enough such that 
\begin{equation}\label{Cond-T0}
 2^{k+1} (c\,\Vert u_0 \Vert_{H^s})^{k} \, C_\eta e^{c_\eta \, T_0}\, T^{1-\alpha}_{0}<1,   
\end{equation} 
and then, the existence and uniqueness of a (local in time) solution $u \in \mathcal{C}([0,T_0], H^s(\R))$ of the integral equation (\ref{Equation-Int}) follow from standard arguments.\finpv
In order to study the regularity (in the spatial variable) of solutions of equation (\ref{Equation}), we define the space $H^{\infty}(\R)$ as $\ds{H^{\infty}(\R)= \bigcap_{r \geq 0}H^r(\R)}$.  
\begin{Proposition}\label{Regular} Let $u \in \mathcal{C}([0,T_0], H^s(\R))$ be the solution of equation (\ref{Equation-Int}) given by Proposition \ref{LocalExistence}. Then, this solution verifies $u \in \mathcal{C}(]0,T_0], H^{\infty}(\R))$.  Moreover we have $u \in \mathcal{C}^{1}(]0,T_0], \mathcal{C}^{\infty}(\R))$ and then $u(t,x)$ is a classical solution of  equation (\ref{Equation}). 
\end{Proposition}	
\pv We will prove that each term in the integral equation (\ref{Equation-Int}) belong to the space $\mathcal{C}(]0,T_0[, H^{\infty}(\R))$. For the first term in the right-hand side of (\ref{Equation-Int}), remark that setting the parameters $s_1 =s$, $s_2>0$ and $\psi = u_0$ in the framework of Lemma \ref{Lemma-tech-1}  then we have  $K_{m,n}(t, \cdot)\ast u_0 \in H^{\infty}(\R)$ pointwise for all $t>0$. Moreover, by Lemma \ref{Lemma-tech-2} we get $K_{m,n}(t, \cdot)\ast u_0 \in \mathcal{C}(]0, T_0],H^{\infty}(\R))$. \\   

We study now the second term in the right-hand side of (\ref{Equation-Int}).  Recall that the solution $u$ of this equation verifies $u(t,\cdot)  \in H^s(\R)$ for all $0 \leq  t \leq T_0$. With this information, and moreover, for $\delta >0$ small enough, first  we will prove  that for all $0 < t \leq T_0$ we have $\ds{ \int_{0}^{t}K_{m,n}(t-\tau,\cdot) \ast u^k \partial_x u (\tau, \cdot) d \tau \in H^{s+\delta}(\R)}$. 
Indeed, in the setting of Lemma \rd{\ref{Lemma-tech-1}}, we set the parameters $s_1=s-1>1/2$, $s_2=\delta+1$  and $\psi= u^k \partial_x u$. Then we write 
\begin{eqnarray*}
& & 	\left\Vert  \int_{0}^{t}K_{m,n}(t-\tau,\cdot) \ast u^k \partial_x u (\tau, \cdot) d \tau \right\Vert_{H^{s+\delta}} \leq \int_{0}^{t} \Vert K_{m,n}(t-\tau,\cdot) \ast u^k \partial_x u (\tau, \cdot) \Vert_{H^{s+\delta}}   \\
&\leq &  \int_{0}^{t} \Vert K_{m,n}(t-\tau,\cdot) \ast u^k \partial_x u (\tau, \cdot) \Vert_{H^{(s-1)+(\delta+1)}} \leq  C_{\eta,\delta}\int_{0}^{t} \frac{ e^{c_{\eta,\delta}(t-\tau)} }{(t-\tau)^{(\delta+1) \alpha}} \Vert u^k \partial_x u (\tau, \cdot) \Vert_{H^{s-1}} d \tau \\
&\leq &  C_{\eta,\delta}e^{c_{\eta,\delta}\,T_0}\int_{0}^{t} \frac{\Vert u^{k}(\tau, \cdot )\Vert_{H^{s-1}} \Vert \partial_x u(\tau, \cdot)\Vert_{H^{s-1}}}{(t-\tau)^{(\delta+1) \alpha}}  d \tau \leq  C_{\eta, \delta, T_0} \int_{0}^{t} \frac{\Vert u(\tau, \cdot )\Vert^{k}_{H^{s-1}} \Vert  u(\tau, \cdot)\Vert_{H^{s}}}{(t-\tau)^{(\delta+1) \alpha}}  d \tau  \\
&\leq &  C_{\eta, \delta, T_0}  \int_{0}^{t} \frac{\Vert u(\tau, \cdot )\Vert^{k+1}_{H^{s}}}{(t-\tau)^{(\delta+1) \alpha}}  d \tau \leq   C_{\eta,\delta, T_0}  \left( \sup_{0 \leq \tau \leq T_0} \Vert u(\tau, \cdot) \Vert_{H^s}\right)^{k+1} \int_{0}^{t} \frac{d \tau}{(t-\tau)^{(\delta+1) \alpha}}.   
\end{eqnarray*}	Here, as $1/\alpha \geq 2$ (see the expression (\ref{Alpha})) then  we set   $0< \delta < 1/\alpha -1$, hence we have $1-(1+\delta)\alpha>0$, and then  the last integral computes down as  $\ds{ \int_{0}^{t} \frac{d \tau}{(t-\tau)^{(\delta+1) \alpha}} \leq c\,  t^{1-(1+\delta)\alpha}}$. Then, for all   $0 < t \leq T_0$ we obtain  
\begin{equation}\label{eq12}
 \left\Vert  \int_{0}^{t}K_{m,n}(t-\tau,\cdot) \ast u^k \partial_x u (\tau, \cdot) d \tau \right\Vert_{H^{s+\delta}}  \leq C_{\eta, \delta, T_0}  \left( \sup_{0 \leq \tau \leq T_0} \Vert u(\tau, \cdot) \Vert_{H^s}\right)^{k+1} \,  t^{1-(1+\delta)\alpha}.
\end{equation}
We prove now  the continuity respect to the temporal variable. Let  $0<t_1, t_2 \leq T_0$ and assume (without loss of generality) that $t_2<t_1$. We write 
\begin{eqnarray*}
& & \left\Vert  \int_{0}^{t_1}K_{m,n}(t_1-\tau,\cdot) \ast u^k \partial_x u (\tau, \cdot) d \tau - \int_{0}^{t_2}K_{m,n}(t_2-\tau,\cdot) \ast u^k \partial_x u (\tau, \cdot) d \tau  \right\Vert_{H^{s+\delta}}\\
&\leq &  \left\Vert  \int_{0}^{t_1}K_{m,n}(t_1-\tau,\cdot) \ast u^k \partial_x u (\tau, \cdot) d \tau - \int_{0}^{t_2}K_{m,n}(t_1-\tau,\cdot) \ast u^k \partial_x u (\tau, \cdot) d \tau  \right\Vert_{H^{s+\delta}} \\
 & & + \left\Vert  \int_{0}^{t_2}K_{m,n}(t_1-\tau,\cdot) \ast u^k \partial_x u (\tau, \cdot) d \tau - \int_{0}^{t_2}K_{m,n}(t_2-\tau,\cdot) \ast u^k \partial_x u (\tau, \cdot) d \tau  \right\Vert_{H^{s+\delta}}\\
 &\leq & \left\Vert \int_{t_2}^{t_1} K_{m,n}(t_1-\tau,\cdot) \ast u^k \partial_x u (\tau, \cdot) d \tau \right\Vert_{H^{s+\delta}} + \left\Vert \int_{0}^{t_2} (K_{m,n}(t_1-\tau, \cdot)-K_{m,n}(t_2-\tau, \cdot))\ast u^k \partial_x u (\tau, \cdot) d \tau \right\Vert_{H^{s+\delta}}.   
\end{eqnarray*}
By (\ref{eq12}) the first term in the right-hand side is estimated as
$$ \left\Vert \int_{t_2}^{t_1} K_{m,n}(t_1-\tau,\cdot) \ast u^k \partial_x u (\tau, \cdot) d \tau \right\Vert_{H^{s+\delta}} \leq C_{\eta, \delta, T_0}  \left( \sup_{0 \leq \tau \leq T_0} \Vert u(\tau, \cdot) \Vert_{H^s}\right)^{k+1} \,  \vert t_1 - t_2 \vert^{1-(1+\delta)\alpha}. $$ 
For the second term in the right-hand side, by Lemma \ref{Lemma-tech-2} we have
\begin{eqnarray*}
& & 	\left\Vert \int_{0}^{t_2} (K_{m,n}(t_1-\tau, \cdot)-K_{m,n}(t_2-\tau, \cdot))\ast u^k \partial_x u (\tau, \cdot) d \tau \right\Vert_{H^{s+\delta}} \leq  c_1  \vert t_1 - t_2 \vert^{1/2} \,  \int_{0}^{t_2}  \Vert u^k \partial_x u (\tau, \cdot) \Vert_{H^{s-1}} d \tau\\
&\leq & c_1  \vert t_1 - t_2 \vert^{1/2} \,  \left( \sup_{0 \leq \tau \leq T_0} \Vert u(\tau, \cdot) \Vert_{H^s}\right)^{k+1}  \, T_0. 
\end{eqnarray*}
 Thus, by these estimates we get $\ds{\int_{0}^{t} K_{m,n}(t-\tau, \cdot)\ast u^k \partial_x u (\tau, \cdot) d \tau \in \mathcal{C}(]0, T_0[, H^{s+\delta}(\R))}$ for   $0< \delta < 1/\alpha -1$. \\
 
 At this point, we have proved that $\ds{u \in  \mathcal{C}(]0, T_0[, H^{s+\delta}(\R))}$ and repeating this process (in order to obtain a gain of regularity for the non linear term) we conclude that $\ds{u \in  \mathcal{C}(]0, T_0[, H^{\infty}(\R))}$. Thereafter, we observe that $u(t,x)$ solves the equation (\ref{Equation}) in the classical way and moreover, writing 
 $$ \partial_t u = - D (\partial_x u) +u^k \partial_x u - \eta(\mathcal{H} \partial^{n}_x u +   \mathcal{H}_m u),$$ we get  $\partial_t u  \in \mathcal{C}(]0, T_0[, H^{\infty}(\R))$. From this information we can  verify now that we have $\partial_{t} u \in \mathcal{C}(]0, T_0[,\mathcal{C}^{\infty}(\R))$. Indeed, we will prove that for all $k\in \mathbb{N}$, the function $\ds{\partial^{k}_{x} \partial_ t u(t,\cdot)}$ is a H\"older continuous function  on $\R$. Let  $k\in \mathbb{N}$ fix. Then, for $\frac{1}{2}<s_1<\frac{3}{2}$ we set   $r=k+s_1$  and since $\partial_t u  \in \mathcal{C}(]0, T_0[, H^{\infty}(\R))$ then  we have $\ds{\partial^{k}_{x} \partial_t u(t,\cdot) \in H^{s_1}(\R)}$. On the other hand, recall that we have the identification  $H^{s_1}(\R)=B^{s_1,2}_{2}(\R)$ (where $B^{s_1,2}_{2}(\R)$ denotes a Besov space \cite{BaChDan}) and moreover we have the inclusion $\ds{B^{s_1,2}_{2}(\R) \subset B^{s_1-\frac{1}{2},\infty}_{\infty}(\R)}$. Thus we get  $\partial^{k}_{x}\partial_t u(t,\cdot) \in \dot{B}^{s_1-\frac{1}{2},\infty}_{\infty}(\R)$.  But, by definition of the space $\ds{\dot{B}^{s_1-\frac{1}{2},\infty}_{\infty}(\R)}$ (see always \cite{BaChDan}) and  since $\frac{1}{2}<s_1<\frac{3}{2}$, then we have $ 0<s_1-\frac{1}{2}<1$ and thus  $\partial^{n}_{x}\partial_t u(t,\cdot)$ is a $\beta$- H\"older continuous function with $\beta =s_1-\frac{1}{2}$. Then we have $\partial_{t} u \in \mathcal{C}(]0, T_0[,\mathcal{C}^{\infty}(\R))$ and thus  $u \in \mathcal{C}^{1}(]0, T_0[,\mathcal{C}^{\infty}(\R))$.\finpv 
Finally, we prove the global well-posedness. Following similar arguments of \bl{\cite{CuiZhao}} (see the proof of Theorem $2$, page $9$) we have the following result. 
\begin{Proposition}\label{Prop-global} Let $T^{*}>0$ be the maximal time of existence of a unique solution $u \in \mathcal{C}([0, T^{*}[, H^s(\R))$ for the equation (\ref{Equation-Int}) given by Proposition \ref{LocalExistence}. Then  we have $T^{*}=+\infty$. 
\end{Proposition}	
\pv By definition we have: $$\ds{T^{*}= \sup \left\{ T>0 : \,\text{there exists a unique solution}\, u\in \mathcal{C}([0,T[, H^s(\R))\,\, \text{of}\,\, (\ref{Equation-Int})\,\,\text{arising from}\,\, u_0 \right\}}.$$  We will assume that $T^{*}<+\infty$ which give us a contradiction. First we need to derive  an energy estimate for solution $u(t,x)$ and for this recall that by  Proposition \rd{\ref{Regular}} we know that this solution is regular enough and then it verify the equation (\ref{Equation})  in a classical way. Thus, we can multiply this equation pointwise by $u$ and integrating in the spatial variable (after some integration by parts) we get:
\begin{equation}\label{eq13}
\frac{1}{2} \frac{d}{d t} \Vert u(t,\cdot) \Vert^{2}_{L^2}= -  \eta\int_{\R} (\mathcal{H} \partial^{n}_x u +   \mathcal{H}_m u) u \, dx,
\end{equation}
  where we must study the term in the right-hand side. By the Parseval's identity we write 
$$- \eta\int_{\R} (\mathcal{H} \partial^{n}_x u +   \mathcal{H}_m u) u \, dx =  - \eta \underbrace{\int_{\R} (i^{n+1} \vert \xi \vert \xi^{n-1} +\vert \xi \vert^m) \vert \widehat{u} \vert^2 d \xi}_{(a)},  $$ and we will estimate the quantity (a) respect to the following values of the parameter $n$ given by (\ref{Cond-n}).  For $n=1$ we have 
\begin{eqnarray*}
(a) &=&   \int_{\R} (\vert \xi \vert - \vert \xi \vert^m) \vert \widehat{u}\vert^2 d\xi = \int_{\vert \xi \vert \leq 2^{1/(m-1)}}  (\vert \xi \vert - \vert \xi \vert^m) \vert \widehat{u}\vert^2 d\xi + \int_{\vert \xi \vert \geq  2^{1/(m-1)}}  (\vert \xi \vert - \vert \xi \vert^m) \vert \widehat{u}\vert^2 d\xi \\
&\leq &  \int_{\vert \xi \vert \leq 2^{1/(m-1)}} \vert \xi \vert  \vert \widehat{u} \vert^2 d \xi -  \int_{\vert \xi \vert \leq 2^{1/(m-1)}} \vert \xi \vert^m  \vert \widehat{u} \vert^2 d \xi  -  \int_{\vert \xi \vert \geq 2^{1/(m-1)}} \vert \xi \vert \vert \widehat{u} \vert^2 d \xi \\
&\leq & c \int_{\vert \xi \vert \leq 2^{1/(m-1)}}  \vert \widehat{u} \vert^2 d \xi \leq c \int_{\R} \vert \widehat{u} \vert^{2} d \xi \leq c \Vert u(t,\cdot) \Vert^{2}_{L^2}. 
\end{eqnarray*} Getting back to (\ref{eq13}) and using the Gr\"owall inequality we have, for all $t \in [0, T^{*}[$, $\ds{\Vert u(t,\cdot) \Vert^{2}_{L^2} \leq c \Vert u_0 \Vert^{2}_{L^2} e^{2 \eta t}}$.\\

Then, for $n=2d$ with $d \in \mathbb{N}^{*}$ we write $\ds{(a)= -\int_{\R} i^{n+1} \vert \xi \vert \xi^{n-1}  \vert \widehat{u} \vert^2 d \xi- \int_{\R} \vert \xi \vert^m \vert \widehat{u} \vert^2 d \xi}$, and since $\ds{i^{n+1} \vert \xi \vert \xi^{n-1}  \vert \widehat{u} \vert^2}$ is a odd function we obtain $\ds{(a)= - \int_{\R} \vert \xi \vert^m \vert \widehat{u} \vert^2 d \xi \leq 0}$. So, by (\ref{eq13}) we have $\ds{\Vert u(t,\cdot) \Vert^{2}_{L^2} \leq c  \Vert u_0 \Vert^{2}_{L^2}\leq  c \Vert u_0 \Vert^{2}_{L^2} e^{2 \eta t}}$, for all $t \in [0, T^{*}[$.\\

Finally, for $n=3+4d$, with $d \in \mathbb{N}$, remark that the function $\ds{i^{n+1} \vert \xi \vert \xi^{n-1} +\vert \xi \vert^m}$ writes down as $\ds{\vert \xi \vert^n +\vert \xi \vert^m }$ and then we have  $\ds{(a)= - \int_{\R} (\vert \xi \vert^n +\vert \xi \vert^m)\vert \widehat{u}\vert^2 d\xi \leq 0}$. Thus, always by (\ref{eq13}) we obtain the same estimate above. \\

We have proven that the solution $u(t,x)$ verifies, for all $t \in [0, T^{*}[$, the energy estimate  
\begin{equation}\label{energ}
 \Vert u(t,\cdot) \Vert_{L^2} \leq c \Vert u_0 \Vert_{L^2} e^{\eta T^{*}},
\end{equation}
 and with this estimate (assuming that $T^{*}<+\infty$)  we will obtain a contradiction as follows: first, we set the fixed quantity $\ds{M= c \Vert u_0 \Vert_{L^2} e^{\eta T^{*}}>0}$. Then,
for  any initial datum $v_0 \in H^{s}(\R)$, recall that  by estimate (\ref{Cond-T0})  the time $T=T(v_0)>0$ of existence of a solution $v \in \mathcal{C}([0,T[, H^s(\R))$ of equation (\ref{Equation-Int}) arising from an initial datum $v_0$ is controlled as: 
$$  2^{k+1} (c\,\Vert u_0 \Vert_{H^s})^{k} \, C_\eta e^{c_\eta \, T}\, T^{1-\alpha}_{0}<1,$$ hence we write 
$$ e^{c_\eta \, T} T^{1-\alpha} < \frac{1}{2^{k+1} C_\eta (c\,\Vert v_0 \Vert_{H^s})^{k}}.$$ 
Moreover, remark that as we have  $\Vert v_0 \Vert_{L^2}\leq  \Vert v_0 \Vert_{H^s}$, then the existence time $T=T(v_0)$ may be controlled by the quantity $\Vert v_0 \Vert_{L^2}$ as follows: 
$$ e^{c_\eta \, T} T^{1-\alpha} < \frac{\eta^\alpha}{2^{k+1} C_\eta (c\,\Vert v_0 \Vert_{L^2})^{k}}.$$ 
In this estimate we may observe that the existence time $T=T (v_0)$ is a \emph{decreasing} function of $\Vert v_0 \Vert_{L^2}$ and then, there exist a time $0<T_1 < T^{*}$ such that for all initial datum $v_0 \in H^{s}(\R)$  such that $\Vert v_0 \Vert_{L^2} \leq M$ the associated solution $v \in  \mathcal{C}([0, T[, H^s(\R))$ exists at least on the interval  $[0,T_1]$ and verifies  $v \in \mathcal{C}([0,T_1], L^2(\R))$.  Thus, for $0< \varepsilon < T_1$ and for the solution $u(t,x)$ (arising from $u_0$)   we can consider the initial datum $v_0= u(T^{*}-\varepsilon, \cdot) \in H^s(\R)$, which by (\ref{energ}) verifies $\Vert v_0 \Vert_{L^2}\leq M$. So, there exists a solution $v$ arising from  $v_0= u(T^{*}-\varepsilon, \cdot)$ which is defined at least on $[0,T_1]$. Thus, gathering  the functions  $u(t,x)$ and $v(t,x)$  we get a solution   
\begin{equation*}
\tilde{u}(t,\cdot) = \left\{ \begin{array}{ll}\vspace{2mm} u(t,\cdot) & \text{for} \,\, 0\leq t \leq T^{*}-\varepsilon, \\
v(t,\cdot) & \text{for}\,\, T^{*}-\varepsilon \leq t \leq T^{*} - \varepsilon +T_1, \end{array}	 \right.
\end{equation*}
arising from the datum $u_0$ which  is defined on the interval $[0, T^{*} - \varepsilon +T_1]$. But, since $0<\varepsilon <T_1$ we have $T^{*} - \varepsilon +T_1 > T^{*}$ and then we contradict the definition of $T^{*}$. Then we have $T^{*}=+\infty$. \finpv 
\section{Pointwise decaying properties}\label{Sec:PWD}
\subsection{Proof of Theorem \rd{\ref{Th-Decay}}}
Given an initial $u_0 \in H^s(\R)$ (with $s>3/2$) by Theorem \rd{\ref{Th-GWP}} there exists a unique solution $u \in \mathcal{C}([0,+\infty[,H^{s}(\R))$ of equation (\ref{Equation}). 
By (\ref{Cond-datum}) we assume now that the initial datum verifies moreover $u_0 \in L^{\infty}((1+\vert \cdot \vert^{\gamma}) dx)$, with $\gamma>0$, and  we will construct a solution $u_1(t,x)$ of equation (\ref{Equation}) in the functional space  $$ \mathcal{C}([0,+\infty[,H^{s}(\R)) \cap E_{\alpha,n, \gamma},$$ where, for  the parameter $0<\alpha \leq 1/2$ defined in (\ref{Alpha}), the parameter $n \geq 1$ given by (\ref{Cond-n})   and the parameter $\gamma>0$ above, the space $E_{\alpha,n, \gamma}$ characterizes the pointwise spatial decaying of solutions  and it is defined as: 
\begin{equation}\label{def-E}
E_{\alpha,\beta,\gamma}= \left\{ u \in \mathcal{S}^{'}([0,+\infty[\times \R):\,\, \text{for all}\,\, 0<T<+\infty,\,\, \sup_{0<t\leq T} t^{\alpha} \Vert (1+\vert \cdot \vert^{\min(\gamma,n+1)}) u(t,\cdot) \Vert_{L^{\infty}} <+\infty \right\}. 
\end{equation}  
It is worth to mention that the weight in the temporal variable $t^\alpha$ is essentially technical (due to the kernel estimates (\ref{Kernel-estimates}))  and it will be useful to carry up all our estimates. We start by the local in time existence of  solutions. More precisely, for a time $0<T<+\infty$, which we shall fix small enough later, we will solve the equation
\begin{equation}\label{Equation-int-aux}
 u_1(t,\cdot)= K_{m,n}(t,\cdot)\ast u_0+ \int_{0}^{t} K_{m,n}(t-\tau) \ast (u^{k}_{1}\partial_x u_{1})(\tau,\cdot) d\tau,   
\end{equation} in the Banach space $\ds{E_T=\mathcal{C}([0,T], H^s(\R))\cap E_{\alpha,n,\gamma}}$, with the norm 
\begin{equation}\label{normET}
  \Vert u \Vert_{E_T}= \sup_{0\leq t\leq T} \Vert u(t,\cdot)\Vert_{H^s}+\sup_{0<t\leq T} t^{\alpha} \Vert (1+\vert \cdot \vert^{\min(\gamma,n+1)})u(t,\cdot) \Vert_{L^{\infty}}.   
\end{equation}

We study the first term in the right-hand side in (\ref{Equation-int-aux}).  Recall that the quantity $\ds{\sup_{0\leq t \leq T} \Vert K_{m,n}(t,\cdot)\ast u_0 \Vert_{H^s}}$ was  estimated in  (\ref{estim-lin1}) and then it remains to  estimate the quantity $\ds{\sup_{0<t\leq T} t^{\alpha} \Vert (1+\vert \cdot \vert^{\min(\gamma,n+1)}) K_{m,n}(t,\cdot)\ast u_0\Vert_{L^{\infty}}}$. As $u_0 \in  L^{\infty}((1+\vert \cdot \vert^{\gamma})dx)$ and moreover, by the kernel estimate (\ref{Kernel-estimates}), for $t>0$ and $x \in \R$ fix we write  
\begin{eqnarray*}
& & \vert K_{m,n}(t,\cdot)\ast u_0(x) \vert  \leq  \int_{\R} \vert K_{m,n}(t, x-y)\vert \vert u_0(y)\vert dy \leq \int_{\R} \vert K_{m,n}(t, x-y)\vert \frac{1+\vert y \vert^{\gamma}}{1+\vert y \vert^{\gamma}} \vert u_0(y)\vert dy\\
&\leq & \Vert (1+ \vert \cdot \vert^{\gamma})  u_0 \Vert_{L^{\infty}} \int_{\R}  \frac{\vert K_{m,n}(t,x-y)\vert }{1+\vert y \vert^{\gamma}} dy \leq \Vert (1+ \vert \cdot \vert^{\gamma})  u_0 \Vert_{L^{\infty}} \,C_\eta  \frac{e^{c_\eta\, t}}{t^{\alpha}} \int_{\R} \frac{dy}{(1+\vert x-y\vert^{n+1})(1+\vert y \vert^{\gamma})} \\
&\leq & \Vert (1+ \vert \cdot \vert^{\gamma})  u_0 \Vert_{L^{\infty}} \,C_\eta \frac{e^{c_\eta\, t}}{ t^{\alpha}} \, \frac{c}{1+\vert x \vert^{\min(\gamma,n+1)}},
\end{eqnarray*}	hence we have 
\begin{equation}\label{Estim-Lin-aux}
t^{\alpha}\left\Vert (1+\vert \cdot \vert^{\min(\gamma,\beta)})K_{m,n}(t,\cdot)\ast u_0 \right\Vert_{L^{\infty}} \leq C_\eta\, e^{c_\eta \, t}  \Vert (1+\vert \cdot\vert^{\gamma}) u_0 \Vert_{L^{\infty}}.  
\end{equation}
Thus, by estimates (\ref{estim-lin1}) and  (\ref{Estim-Lin-aux})  we get
\begin{equation}\label{Estim-Lin}
\Vert K_{m,n}(t,\cdot) \ast u_0 \Vert_{E_T} \leq C_\eta\,  e^{c_\eta\, T} ( \Vert (1+\vert \cdot \vert^{\gamma} )  u_0 \Vert_{L^{\infty}} + \Vert u_0 \Vert_{H^s}).
\end{equation}
We study now the second term in the right-hand side in the equation (\ref{Equation-int-aux}). Remark that  since the quantity $\ds{\sup_{0<t\leq T} \left \Vert  \int_{0}^{t} K_{m,n}(t-\tau, \cdot)\ast (u^{k}_{1} \partial_x u_1)(\tau, \cdot) d \tau \right\Vert_{H^s}}$ was estimated in (\ref{estim-nonlin1}) it remains to estimate the quantity $\ds { \sup_{0<t\leq T} t^{a} \left \Vert (1+\vert \cdot\vert^{\min(\gamma,n+1)}) \left( \int_{0}^{t} K_{m,n}(t-\tau, \cdot)\ast (u^k \partial_x u)(\tau, \cdot) d \tau  \right)\right\Vert_{L^{\infty}}}$. By the kernel estimate (\ref{Kernel-estimates}), for $0<\tau <t \leq T$ and $x \in \R$ fix we write: 
\begin{eqnarray*} 
& & \vert K_{m,n}(t-\tau, \cdot)\ast (u^{k}_{1} \partial_x u_1)(\tau, x) \vert \leq \int_{\R} \vert K_{m,n}(t-\tau, x-y)\vert \vert u_1(\tau, y)\vert^k \vert \partial_y u_1 (\tau, y) \vert dy \\ \nonumber
&\leq &C_\eta \int_{\R} \frac{e^{c_\eta (t-\tau)}}{ (t-\tau)^{\alpha}} \frac{1}{(1+ \vert x-y \vert^{n+1})} \,   	\vert u_1(\tau, y)\vert^k \vert \partial_y u_1 (\tau, y) \vert dy \\ \nonumber
&\leq & C_\eta  \frac{e^{c_\eta (t-\tau)}}{ (t-\tau)^{\alpha}}  \Vert u_1(\tau, \cdot) \Vert^{k-1}_{L^{\infty}} \Vert \partial_x u_1 (\tau, \cdot) \Vert_{L^{\infty}} \, \int_{\R}  \frac{1}{(1+ \vert x-y \vert^{n+1})}  \vert u_1(\tau, y ) \vert dy \\ \nonumber 
&\leq & C_\eta \frac{e^{c_\eta (t-\tau)}}{(t-\tau)^{\alpha}}  \Vert u_1(\tau, \cdot) \Vert^{k-1}_{L^{\infty}} \Vert \partial_x u_1 (\tau, \cdot) \Vert_{L^{\infty}} \Vert (1+\vert \cdot \vert^{\min(\gamma,n+1)}) u_1(\tau, \cdot)\Vert_{L^{\infty}}\\ \nonumber
& & \times  \int_{\R}  \frac{1}{(1+ \vert x-y \vert^{n+1})(1+\vert y \vert^{\min(\gamma,n+1)})} dy \\
&\leq & C_\eta \frac{e^{c_\eta (t-\tau)}}{(t-\tau)^{\alpha}} \underbrace{\Vert u_1(\tau, \cdot) \Vert^{k-1}_{L^{\infty}} \Vert \partial_x u_1 (\tau, \cdot) \Vert_{L^{\infty}} \Vert (1+\vert \cdot \vert^{\min(\gamma,n+1)}) u_1(\tau, \cdot)\Vert_{L^{\infty}} }_{(a)}  \, \frac{1}{1+\vert x \vert^{\min(\gamma,n+1)}}. 
\end{eqnarray*} In the last expression we still  must estimate the term  $(a)$. Recall that as $s>3/2$ then we have  $s-1>1/2$ and thus the space $H^{s-1}(\R)$  embeds in the space $L^{\infty}(\R)$. So we can write  
\begin{eqnarray}\label{estim_a} \nonumber 
(a)& \leq & \Vert u_1(\tau, \cdot) \Vert^{k-1}_{H^{s-1}} \Vert \partial_x u_1(\tau, \cdot) \Vert_{H^{s-1}} \frac{\tau^{\alpha}}{\tau^{\alpha}} 	\Vert (1+\vert \cdot \vert^{\min(\gamma,n+1)}) u_1(\tau, \cdot)\Vert_{L^{\infty}}\\
&\leq & \frac{c}{\tau^{\alpha}} \Vert u(\tau, \cdot) \Vert^{k}_{H^s} \left( \tau^{\alpha} \Vert(1+\vert \cdot \vert^{\min(\gamma,n+1)}) u_1(\tau, \cdot)\Vert_{L^{\infty}}\right).
\end{eqnarray}
Thus, getting back to the previous estimate we have 
\begin{equation*}
\begin{split}
\Vert (1+\vert x \vert^{\min(\gamma,n+1)} )\vert K_{m,n}(t-\tau, \cdot)\ast (u^{k}_{1} \partial_x u_1)(\tau, \cdot) \Vert_{L^{\infty}} \leq & \,  C_\eta \frac{e^{c_\eta (t-\tau)}}{(t-\tau)^{\alpha} \tau^{\alpha}} \Vert u(\tau, \cdot) \Vert^{k}_{H^s}\\
& \times \left( \tau^{\alpha} \Vert(1+\vert \cdot \vert^{\min(\gamma,n+1)}) u_1(\tau, \cdot)\Vert_{L^{\infty}}\right),    
\end{split}
\end{equation*}  
hence we can write 
\begin{equation}\label{I2gronwall}
\begin{split}
&  t^{\alpha} \int_{0}^{t} \left \Vert (1+\vert \cdot\vert^{\min(\gamma,n+1)}) \left(  K_{m,n}(t-\tau, \cdot)\ast (u^{k}_{1} \partial_x u_1)(\tau, \cdot) \right)\right\Vert_{L^{\infty}}  d \tau  \\
\leq & C_\eta\, e^{c_\eta\, T}   t^{\alpha} \int_{0}^{t} \frac{1}{(t-\tau)^\alpha \, \tau^{\alpha}}  \Vert u(\tau, \cdot) \Vert^{k}_{H^s} \left( \tau^{\alpha} \Vert(1+\vert \cdot \vert^{\min(\gamma,n+1)}) u_1(\tau, \cdot)\Vert_{L^{\infty}}\right) d \tau.
\end{split}
\end{equation}
Now, recalling the definition of the norm $\Vert \cdot \Vert_{E_T}$ given in (\ref{normET}) we finally get 
$$ \sup_{0<t\leq T} t^{\alpha}\left\Vert (1+\vert \cdot\vert^{\min(\gamma,n+1)}) \int_{0}^{t} K_{m,n}(t-\tau, \cdot)\ast (u^{k}_{1} \partial_x u_1)(\tau, \cdot) d \tau \right\Vert_{L^{\infty}}  \leq C_\eta\, e^{c_\eta\, T} T^{1-\alpha} \Vert u_1 \Vert^{k+1}_{E_T}.$$
By this estimate and by estimate (\ref{estim-nonlin1}) we obtain
\begin{equation}\label{eq10}
\left \Vert \int_{0}^{t} K_{m,n}(t-\tau, \cdot)\ast (u^{k}_{1} \partial_x u_1)(\tau, \cdot) d \tau \right\Vert_{E_T} \leq C_\eta \, e^{c_\eta \, T} \, T^{1-\alpha} \,  \Vert u_1 \Vert^{k+1}_{E_T}. 
\end{equation}
Once we have the estimates (\ref{Estim-Lin}) and (\ref{eq10}), for a time $0<T_0<+\infty$ small enough,  the existence and uniqueness of a (local in time) solution $u_1 \in E_T$ of the integral equation (\ref{Equation-int-aux}) follow from well-known arguments.\\ 

Now  we will show that the solution $u_1(t,x)$ is global in time. Recall that by Proposition \rd{\ref{Prop-global}}  we have $u_1 \in \mathcal{C}([0,+\infty[,H^s(\R))$ and then it remains to prove that the quantity $\ds{\sup_{0<t\leq T} t^\alpha \Vert (1+\vert \cdot \vert^{\min(\gamma,n+1)}) u_1(t,\cdot) \Vert_{L^{\infty}}}$ is well-defined for all time $T>0$. \\

Let $T>0$. For all $0<t \leq T$, let us define the quantity $\ds{g(t)=t^\alpha \Vert (1+\vert \cdot \vert^{\min(\gamma,n+1)}) u_1(t,\cdot) \Vert_{L^{\infty}}}$,  and  by equation  (\ref{Equation-int-aux}) we write 
\begin{eqnarray*}
 g(t)&=& t^{\alpha}\left\Vert (1+\vert \cdot\vert^{\min(\gamma,n+1)}) \left( K_{m,n}(t,\cdot)\ast u_0+ \int_{0}^{t} K_{m,n}(t-\tau) \ast (u^{k}_{1}\partial_x u_{1})(\tau,\cdot) d\tau \right) \right\Vert_{L^{\infty}}\\
 &\leq & t^{\alpha}\left\Vert (1+\vert \cdot\vert^{\min(\gamma,n+1)}) K_{m,n}(t,\cdot)\ast u_0 \right\Vert_{L^{\infty}}\\
 & & + t^\alpha \int_{0}^{t} \left\Vert (1+\vert \cdot\vert^{\min(\gamma,n+1)}) \left( K_{m,n}(t-\tau) \ast (u^{k}_{1}\partial_x u_{1})(\tau,\cdot) d\tau \right) \right\Vert_{L^{\infty}}= I_1+ I_2,
\end{eqnarray*}
where we must estimate the terms $I_1$ and $I_2$. For $I_1$, by estimate (\ref{Estim-Lin-aux}) we have directly the estimate 
\begin{equation}\label{estim-gron1}
I_1 \leq C_\eta e^{c_\eta \, T} \Vert (1+\vert \cdot\vert^{\gamma}) u_0 \Vert_{L^{\infty}}=C_1(u_0,T).    
\end{equation} For $I_2$, by estimate (\ref{I2gronwall})  and recalling the definition of the expression $g(t)$ given above, we have
\begin{equation}\label{estim-gron2} 
\begin{split}
I_2 & \leq C_\eta e^{c_\eta\,T}\,  \left( \sup_{0\leq \tau \leq t}  \Vert u_1 (\tau, \cdot)\Vert^{k}_{H^s}\right)   \,t^{\alpha}\, \int_{0}^{t} \frac{1}{(t-\tau)^{\alpha} \tau^{\alpha}}g(\tau) d \tau\\
&\leq   C_\eta e^{c_\eta\,T}\,   \left( \sup_{0\leq \tau \leq T}  \Vert u_1 (\tau, \cdot)\Vert^{k}_{H^s}\right)  t^{\alpha}\,    \int_{0}^{t} \frac{1}{(t-\tau)^{\alpha} \tau^{\alpha}}g(\tau) d \tau \\
&=C_2(u,T) \, t^{\alpha}\,    \int_{0}^{t} \frac{1}{(t-\tau)^{\alpha} \tau^{\alpha}}g(\tau) d \tau.
\end{split}
\end{equation}
At this point, we need to distinguish two cases for the parameter  $0<\alpha \leq 1/2$. 
\begin{enumerate}
\item[$\bullet$]\textbf{For $0<\alpha < 1/2$}. By estimates (\ref{estim-gron1}) and (\ref{estim-gron2}),  for all $t\in ]0,T]$  we obtain the following inequality: $$g(t)\leq C_1(u_0,T)+ C_2(u,T) T^\alpha\, \int_{0}^{t} \frac{1}{(t-\tau)^{\alpha}\,\tau^{\alpha}}g(\tau).$$  Now, in order to get  control (global in time) on the quantity $g(t)$ we will use the following technical result.  For a proof  see the Lemma $7.1.2$ of the book \bl{\cite{DHenry}}.
\begin{Lemme}[Gr\"onwall's type inequality I]\label{Lemme-growall-tech} Let $\mathfrak{a}>0$ and $\mathfrak{b}>0$, such that $\mathfrak{a} + \mathfrak{b}>1$. Let  $g:[0,T]\longrightarrow [0,+\infty[$ be a  function such that verifies:
\begin{enumerate}
\item[a)] $g \in L^{1}_{loc}([0,T])$,  $t^{\mathfrak{b}-1} g \in L^{1}_{loc}([0,T])$, and 
\item[b)] there exists two constants $\mathfrak{C}_1> 0$ and $\mathfrak{C}_2>0$,  such that for almost all  $t\in [0,T]$, we have \begin{equation*}
g(t)\leq \mathfrak{C}_1 +\mathfrak{C}_2 \int_{0}^{t}(t-\tau)^{\mathfrak{a}-1}\tau^{\mathfrak{b}-1}g(\tau) d\tau. 
\end{equation*}
\end{enumerate}	
Then, the following statements hold:   
\begin{enumerate}
\item[$1)$] There exists a continuous and increasing function $\varTheta:[0,+\infty[ \longrightarrow [0,+\infty[$, defined by \begin{equation}\label{omeg} 
\ds{\varTheta(t)= \sum_{k=0}^{+\infty}c_k\, t^{\sigma k}},
			\end{equation}
			where $\sigma=\mathfrak{a}+\mathfrak{b}-1>0$ and  moreover, for the  Gamma function $\Gamma (\cdot)$ the coefficients $c_k>0$  are given by the recurrence formula:  $\ds{c_0=1}$ and $\ds{\frac{c_{k+1}}{c_k}= \frac{\Gamma(k \sigma +1)}{\Gamma(k\sigma +\mathfrak{a} +\mathfrak{b})}}$ for $k\geq 1$.
			\item[$2)$] For all time $t\in [0,T]$, we have  $\ds{g(t) \leq c_1 \varTheta \left( c_2^{\frac{1}{\sigma}} \,t \right)}$.\\
 		\end{enumerate}		
	\end{Lemme}
In  this lemma we set the parameters $\mathfrak{a}=1-\alpha>0$ and $\mathfrak{b}=1-\alpha>0$, hence, as $0<\alpha < 1/2$ then we have $\mathfrak{a}+\mathfrak{b}>1$. Moreover,  it is easy to see that  points $a)$ and $b)$ above  are verified, where, in point $b)$ we set the constants $\ds{\mathfrak{C}_1= C_1(u_0,T)}$ and $\ds{\mathfrak{C}_2= C_2(u,T)T^\alpha}$. Thus, by point $2)$  for $\sigma=1-2\alpha>0$ and for all time $t\in ]0,T]$ we obtain  the control  $\ds{ g(t)=  t^\alpha \Vert (1+\vert \cdot \vert^{\min(\gamma,n+1)}) u_1(t,\cdot) \Vert_{L^{\infty}} \leq   c_1\, \varTheta \left( c_2^{\frac{1}{1-\alpha}}\, t \right)}$, hence the quantity $g(t)$ does not explode in a finite time. Thus we have  $u_1 \in \mathcal{C}([0,+\infty[,H^{s}(\R)) \cap E_{\alpha,n, \gamma}$. 
\item[$\bullet$] \textbf{For $\alpha = 1/2$.}  Observe that   by estimates (\ref{estim-gron1}) and (\ref{estim-gron2}),  for all $0<t\leq T$  we obtain the following inequality:
 \begin{equation}\label{estim-gron-critica}
 g(t)\leq C_1(u_0,T)+ C_2(u,T) t^\alpha\, \int_{0}^{t} \frac{1}{(t-\tau)^{1/2}\,\tau^{1/2}}g(\tau).
 \end{equation}
However, we may observe that this case ir more delicate since if in  Lemma  \ref{Lemme-growall-tech} we set $\mathfrak{a}=\mathfrak{b}=1/2$, then the required condition $\mathfrak{a}+\mathfrak{b}>1$ is not verified.  To contour this  problem, we shall use here another argument.

We suppose that $T^{*<+\infty}$ and we will obtain a contradiction.

version of a Gr\"onwall's type inequality. For a proof of this technical result see the  Lemma $3.4$ in \cite{Brandolese}.
\begin{Lemme}[Gr\"onwall's type inequality II]\label{Lemme-growall-tech-2}  Let $g:[0,T]\longrightarrow [0,+\infty[$ be a non-negative and locally bounded function such that, for all $t \in ]0,T]$ it verifies:  
\begin{equation*}
g(t)\leq \mathfrak{C}_1 + \mathfrak{C}_2 \int_{0}^{t} \frac{1}{(t-\tau)^{1/2} \tau^{1/2}} g(\tau) d\tau, 
\end{equation*} for two constants $\mathfrak{C}_1, \mathfrak{C}_2>0$ depending on $T$. Moreover, let $\ds{K= \int_{0}^{1}\frac{d\, \tau}{(1-\tau)^{1/2} \tau^{1/2}}}$. If $\mathfrak{C}_2< 1/K$,  then for all $t\in ]0,T]$ we have  $\ds{g(t)\leq \mathfrak{C}_1}$. 
\end{Lemme}	
In this lemma, we set the constant $\mathfrak{C}_1=C_1(u_0,T)$. Moreover, in the second term of the right in (\ref{estim-gron-critica}) we remark that  we can set a time $0<T_1<T/2$ such that
$\ds{\mathfrak{C}_2=C_2(u,T) T^{\alpha}_{1}<1/K}$. Thus, for all $t\in ]0,T_1]$  we have $\ds{g(t) \leq  \mathfrak{C}_1}$. Thereafter,  we observe  that we can repeat this process as follows: we consider now the initial datum $u(T_1,\cdot)$ and then for all time $t \in [0,T_1]$ we denote by $\tilde{u}_1(t,x)$ the solution of the problem 
$$ \tilde{u}_1(t,x)= K_{m,n}(t,\cdot)\ast u_1(T_1, x)+ \int_{0}^{t} K_{m,n}(t-\tau)\ast (\tilde{u}^{k}_{1}\partial_x \tilde{u}_1)(\tau, x) d\, \tau.$$
Then, for the quantity  $\tilde{g}(t)=t^{\alpha} \Vert (1+\vert \cdot\vert^{\min(\gamma, n+1)})\tilde{u}_1(t,\cdot) \Vert_{L^{\infty}}$, by estimates (\ref{estim-gron1}) and (\ref{estim-gron2}), and moreover, applying the Lemma \ref{Lemme-growall-tech-2}, for all $t\in [0,T_1]$ we have $\ds{\tilde{g}(t)\leq C_1(u(T_1,\cdot),T)}$. But, by uniqueness of solutions in the space $\mathcal{C}([0,+\infty[,H^s(\R))$ we have the identity $\tilde{u}_1(t,x)=u_1(t+T_1,x)$ and then we have $\ds{\tilde{g}(t)=g(t+T_1)\leq C_1(u(T_1,\cdot),T)}$. Repeating this process a finite number of iterations  we arrive to the time  $T>0$. Then (when $\alpha=1/2$) the  quantity $g(t)$ does not explode in a finite time and thus  have $u_1 \in \mathcal{C}([0,+\infty[,H^{s}(\R)) \cap E_{\alpha,n, \gamma}$. 
\end{enumerate}

In order to finish this proof, remark that always by uniqueness of solution $u(t,x)$ in the space $\mathcal{C}([0,+\infty[,H^{s}(\R))$ we have the identity $u=u_1$ and then the solution $u(t,x)$ belongs to the space $E_{\alpha,n,\gamma}$. By definition of the space $E_{\alpha,n,\gamma}$ given in (\ref{def-E}), for all $t>0$ and for $x\in \R$ we can write 
\begin{equation}\label{Decay1-pres}
\vert u(t,x) \vert \leq \frac{1}{t^{\alpha}} \left(\sup_{0<\tau \leq t} \tau^{\alpha} \Vert (1+\vert \cdot\vert^{\min(\gamma, n+1)}) u(\tau, \cdot) \Vert_{L^{\infty}} \right)\frac{1}{1+\vert x \vert^{\min(\gamma,n+1)}} =  \frac{\mathfrak{C}_0(u)}{t^{\alpha}} \frac{1}{1+\vert x \vert^{\min(\gamma,n+1)}},  
\end{equation}
hence, setting the constant $\ds{C_0(u,t)= \frac{\mathfrak{C}_0(u)}{t^{\alpha}}>0}$,   we obtain the desired estimate (\ref{Decay1}).  Theorem \ref{Th-Decay} is proven. \finpv 
\subsection{Proof of Theorem \rd{\ref{Th-Decay2}}} 
Recall that for $n\geq 1$  given by (\ref{Cond-n})  for $\varepsilon \in ]0,1]$ we have $\gamma= n+1+\varepsilon$.  Since   that the solution $u(t,x)$ writes down as in the integral formulation (\ref{Equation-Int}), we start by proving that the first term in the right-hand side in   (\ref{Equation-Int} has the  the following asymptotic development: 
\begin{equation}\label{eq18}
K_{m,n}(t,\cdot)\ast u_0 (x)= K_{m,n}(t,x)\left(  \int_{\R} u_0 (y)dy\right)+ R_{1}(t,x),  \quad \vert x \vert \to +\infty,   
\end{equation} with $\ds{\vert R_1(t,x) \vert \leq  \frac{c(u_0,t)}{\vert x \vert^{\gamma}}}$. Indeed,  for $t>0$ and $x \in \R$ fix  this term can be decomposed as follows: 
\begin{eqnarray*}
\int_{\R} K_{m,n}(t,x-y)u_0(y) dy &=& K_{m,n}(t,x)\left(  \int_{\R} u_0 (y)dy\right)  + \int_{\vert y \vert < \frac{\vert x \vert}{2}} (K_{m,n} (t, x-y)- K_{m,n} (t, x)) u_0(y) dy  \\
& & + \int_{\vert y \vert > \frac{\vert x \vert}{2}} K_{m,n}(t, x-y) u_0(y) dy - K_{m,n}(t,x) \left( \int_{\vert y \vert > \frac{\vert x \vert}{2}} u_0 (y) dy\right)\\
&=&  K_{m,n}(t,x)\left(  \int_{\R} u_0 (y)dy\right)  + I_1 + I_2 + I_3,
\end{eqnarray*}
hence, we define $\ds{R_1=  I_1 + I_2 + I_3 }$ and we will verify that the following statement holds:
\begin{equation}\label{eq14}
\vert R_1 \vert \leq \frac{c(u_0,t)}{\vert x \vert^{\gamma}}, \quad \vert x \vert \to +\infty. 
\end{equation} 
To estimate the term $I_1$ we need the following technical result. Its proof follows the same lines of the proof of Lemma \rd{$4.2$} in \bl{\cite{CorJa}}.  
\begin{Lemme}\label{Lema-der-nucleo} Let the parameter $0<\alpha \leq 1/2$ given in (\ref{Alpha}). Within the framework of Proposition  \ref{Prop-Kernel},  for all time $t>0$ this kernel $K_{m,n}(t,x)$ satisfies $K_{m,n}(t,\cdot) \in \mathcal{C}^{1}(\R)$ and we have the following estimates: 
\begin{enumerate}
\item[$1)$] For all  $x \neq 0$, $\ds{\vert \partial_x K_{m,n}(t,x) \vert \leq C_\eta \frac{e^{c_\eta\,  t }}{\vert x \vert^{n +2}}}$. 
\item[$2)$] For all $x \in \R$, $\ds{\vert \partial_x K_{m,n} (t,x) \vert \leq C_\eta \frac{e^{c_\eta\, t}}{t^{2 \alpha}}\frac{1}{1+ \vert x \vert^{n+2}}}$. 
\end{enumerate}		
\end{Lemme}	
 As  $K_{m,n}(t,\cdot) \in \mathcal{C}^{1}(\R)$,   by the  Taylor expansion we write $ \ds {K_{m,n}(t, x-y) - K_{m,n} (t,x) = - y \partial_x K_{m,n} (t, x- \theta y )}$ for some  $0< \theta < 1$. Then, by this identity and the estimate given in point $1)$ above,   in the term $I_1$ we obtain 
 \begin{eqnarray*}
 I_1 & \leq & 	\int_{\vert y \vert < \frac{\vert x \vert}{2}} \vert  (K_{m,n} (t, x-y)- K_{m,n} (t, x)) \vert \vert  u_0(y) \vert  dy \leq \int_{\vert y \vert < \frac{\vert x \vert}{2}} \vert y \vert \vert \partial_x K_{m,n} (t, x- \theta y ) \vert \vert u_0(y) \vert d y \\
 & \leq & C_\eta e^{c_\eta \, t } \int_{\vert y \vert < \frac{\vert x \vert}{2}} \frac{\vert y \vert \vert u_0 (y) \vert}{\vert x - \theta y \vert^{n+2}} dy.  
 \end{eqnarray*}	
 We study now the expression $\ds{\frac{1}{\vert x -\theta y \vert^{\beta +1}}}$. As we have $0<\theta < 1$ and moreover, as we have $\vert y \vert < \frac{\vert x \vert }{2}$, then we can write $\vert x -\theta y \vert \geq \vert x \vert  - \theta \vert y \vert \geq \vert x \vert - \vert y \vert \geq \frac{\vert x \vert}{2}$; and thus  we get $\ds{\frac{1}{\vert x -\theta y \vert^{n+2}} \leq c \frac{1}{\vert x \vert^{n+2}}}$. With this inequality and recalling that the initial datum verifies $\ds{\vert u_0 (y) \vert \leq \frac{c}{1+\vert y \vert^{\gamma }}}$ (with $\gamma=n+1+\varepsilon$) we can write 
 $$ C_\eta  e^{c_\eta \,  t } \int_{\vert y \vert < \frac{\vert x \vert}{2}} \frac{\vert y \vert \vert u_0 (y) \vert}{\vert x - \theta y \vert^{n +2}} dy \leq \frac{C_\eta e^{c_\eta \,  t }}{\vert x \vert^{n + 2}} \int_{\vert y \vert < \frac{\vert x \vert}{2}} \frac{\vert y \vert }{1+ \vert y \vert^{n+1 +\varepsilon}} dy  \leq \frac{C_\eta e^{c_\eta \, t }}{\vert x \vert^{n + 2}} \int_{\R} \frac{\vert y \vert }{1+ \vert y \vert^{n+1+ \varepsilon }} dy \leq C_\eta \frac{ e^{c_\eta \, t }}{\vert x \vert^{n+ 2}}.$$
Thus, as $\gamma \leq n +2$ then   we have
 \begin{equation}\label{eq15}
 I_1 \leq \frac{C_\eta e^{c_\eta \, t }}{\vert x \vert^{n + 2}}\leq \frac{C_\eta e^{c_\eta \, t }}{\vert x \vert^{\gamma}}, \quad \vert x \vert \to +\infty.   
\end{equation}
 For the term $I_2$, as we have $\ds{\vert u_0 (y) \vert \leq \frac{c}{\vert y \vert^{\gamma}}}$ (for $\vert y \vert$ large enough) and moreover, as we have $\vert y \vert > \frac{\vert x \vert}{2}$, then  we write 
\begin{eqnarray*}
 I_2 &\leq & \int_{\vert y \vert > \frac{\vert x \vert}{2}} \vert K_{m,n}(t,x-y) \vert \vert u_0 (y) \vert d y \leq  c \int_{\vert y \vert > \frac{\vert x \vert}{2}} \frac{\vert K_{m,n}(t, x-y) \vert}{\vert y \vert^{\gamma}} dy \leq \frac{c}{\vert x \vert^{\gamma}} \int_{\vert y \vert > \frac{\vert x \vert}{2}} \vert K_{m,n}(t, x-y) \vert \\
 &\leq & \frac{c}{\vert x \vert^{\gamma}} \Vert K_{m,n} (t,\cdot) \Vert_{L^1},
\end{eqnarray*} but,  by (\ref{Kernel-estimates}) we have $\Vert K_{m,n} (t,\cdot) \Vert_{L^1} \leq C_\eta \frac{e^{c_\eta \, t }}{ t^{\alpha}}$, and then  we can write 
\begin{equation}\label{eq16}
I_2 \leq  C_\eta \frac{e^{c_\eta \,t } }{t^\alpha} \frac{1}{\vert x \vert^{\gamma}}, \quad \vert x \vert \to +\infty.   
\end{equation}	 
Finally, in order to study  the term $I_3$,  recall first  that always  by  (\ref{Kernel-estimates})  for $\vert x \vert $ enough enough we have $\ds{\vert K_{m,n} (t,x) \vert \leq  C_\eta  \frac{e^{c_\eta\, t }}{t^{\alpha}} \frac{1}{\vert x \vert^{n+1}}}$. Moreover, recall that the initial datum verifies $\ds{\vert u_0 (y) \vert \leq  \frac{c}{1+ \vert y \vert^{\gamma}}}$ (with $\gamma=n+1+ \varepsilon$). Then we write      
\begin{eqnarray}\label{eq17}\nonumber
I_3 &\leq &  C_\eta \frac{e^{c_\eta \, t }}{ t^{\alpha}} \frac{1}{\vert x \vert^{n+1}} \int_{\vert y \vert > \frac{\vert x \vert}{2}} \vert u_0(y) \vert dy  \leq  C_\eta\frac{e^{c_\eta \, t }}{ t^{\alpha}} \frac{1}{\vert x \vert^{n+1}} \int_{\vert y \vert > \frac{\vert x \vert}{2}} \frac{1}{1+\vert y \vert^{n+1+\varepsilon}} dy \\
& \leq & C_\eta \frac{e^{c_\eta \, t }}{t^{\alpha}} \frac{1}{\vert x \vert^{n+1+\varepsilon}} \int_{\vert y \vert > \frac{\vert x \vert}{2}} \frac{1}{1+\vert y \vert^{n+1}} dy \leq  C_\eta \frac{e^{c_\eta \, t }}{ t^{\alpha}} \frac{1}{\vert x \vert^{\gamma}} \int_{\R} \frac{1}{1+\vert y \vert^{n+1}} dy \leq C_\eta \frac{e^{c_\eta \, t }}{ t^{\alpha}} \frac{1}{\vert x \vert^{\gamma}}. 
\end{eqnarray}
Thus, the desired estimate (\ref{eq14}) follows from (\ref{eq15}), (\ref{eq16}) and (\ref{eq17}), and we have the asymptotic development given in  (\ref{eq18}). \\

We prove now that the second term in the right-hand side of  (\ref{Equation-Int}) verifies: 
\begin{equation}\label{eq19}
\left\vert \int_{0}^{t} K_{m,n}(t-\tau,\cdot) \ast (u^{k}\partial_x u)(\tau, x) d\tau \right\vert  \leq \frac{\mathfrak{C}_1(u,t)}{\vert x \vert^{n+2}}, \quad \vert x \vert \to +\infty.   
\end{equation} 
For $t>0$ and $x \in \R$ fix we write  
\begin{eqnarray*}
& &  \int_{0}^{t} K_{m,n}(t-\tau) \ast (u^{k}\partial_x u)(\tau, x) d\tau =   \frac{1}{k+1} \int_{0}^{t}  \partial_x K_{m,n}(t-\tau, \cdot) \ast  u^{k+1}(\tau, x) d \tau \\
&=&  \frac{1}{k+1} \int_{0}^{t} \int_{\R} \partial_x K_{m,n}(t-\tau, x-y) u^{k+1}(\tau, y) dy \\
&= &\frac{1}{k+1} \int_{0}^{t} \int_{\R} \partial_x K_{m,n}(t-\tau, x-y) u^2(\tau, y ) u^{k-1}(\tau, y) dy=(a). 
\end{eqnarray*} Then,  by point $2)$ of Lemma \rd{\ref{Lema-der-nucleo}} and recalling that by  estimate (\ref{Decay1-pres}) (with $\gamma=n+1+\varepsilon$) we have the pointwise estimate: $\ds{\vert u(\tau, y) \vert^2 \leq  \frac{\mathfrak{C}^{2}_{0}(u)}{\tau^ {2\alpha} \, (1+\vert y \vert^{2(n+1)})}}$, we obtain 
\begin{eqnarray*}
(a)&\leq & C_\eta \frac{\mathfrak{C}^{2}_{0}(u)}{k+1} \int_{0}^{t} \frac{e^{c_\eta (t-\tau) }}{(t-\tau)^{2 \alpha} \, \tau^{2 \alpha}} \int_{\R} \frac{1}{1+ \vert x-y \vert^{n+ 2}} \frac{1}{1+\vert y \vert^{2(n+1)}}   \vert u^{k-1} (\tau, y ) \vert d y \, d \tau \\
& \leq &  C_\eta \frac{\mathfrak{C}^{2}_{0}(u)}{k+1} \int_{0}^{t} \frac{e^{c_\eta (t-\tau) }}{(t-\tau)^{2 \alpha} \, \tau^{2 \alpha}}  \Vert u (\tau, \cdot) \Vert^{k-1}_{L^{\infty}}  d \tau  \left( \int_{\R} \frac{1}{1+ \vert x-y \vert^{n + 2}} \frac{1}{1+\vert y \vert^{2(n+1) }}   d y  \right)  \\
&\leq & C_\eta \frac{\mathfrak{C}^{2}_{0}(u)\, e^{c_\eta\, t}}{k+1}  \int_{0}^{t} \frac{  \Vert u (\tau, \cdot) \Vert^{k-1}_{L^{\infty}} }{(t-\tau)^{2 \alpha} \, \tau^{2 \alpha}} d \tau  \left( \int_{\R} \frac{1}{1+ \vert x-y \vert^{n + 2}} \frac{1}{1+\vert y \vert^{2(n+1) }}  d y  \right)=(b).
\end{eqnarray*}
Now, recall that as $s > 3/2$ then $H^s(\R)$ embeds in $L^{\infty}(\R)$ and then we can write 
\begin{eqnarray*}
(b)& \leq  &  C_\eta \frac{\mathfrak{C}^{2}_{0}(u)\, e^{c_\eta\, t}}{k+1}   \int_{0}^{t} \frac{  \Vert u (\tau, \cdot) \Vert^{k-1}_{H^s} }{(t-\tau)^{2 \alpha} \, \tau^{2 \alpha}} d \tau  \left( \frac{1}{1+\vert x \vert^{n +2}} \right)  \\
&\leq & C_\eta \frac{\mathfrak{C}^{2}_{0}(u)\, e^{c_\eta\, t}}{k+1} \left(  \sup_{0<\tau < t}  \Vert u(\tau, \cdot) \Vert_{H^s}\right)^{k-1} \left(  \int_{0}^{t} \frac{ d \tau  }{(t-\tau)^{2 \alpha} \, \tau^{2 \alpha}}  \right) \, \frac{1}{1+\vert x \vert^{n +2}}. 
\end{eqnarray*} 
At this point we must estimate the integral in the temporal variable. For this recall the assumption on the parameters $m$ and $n$:  $(m,k)\neq (2,1)$ and $(m,n)\neq (2, 2 d)$ with $d \in \mathbb{N}^{*}$. For those values of $m$ and $n$,  by definition of parameter $0<\alpha\leq 1/2$ given in (\ref{Alpha}) we obtain $0 < \alpha \leq 1/3$ and then this integral computes down as $\ds{  \int_{0}^{t} \frac{ d \tau  }{(t-\tau)^{2 \alpha} \, \tau^{2 \alpha}}  \leq \frac{c}{t^{4 \alpha -1}}}$. At this point, before to continue with the proof of this theorem it is worth to do the following remark. 
\begin{Remarque}\label{CondParam} For the values  $(m,n)= (2,1)$ or $(m,n) = (2, 2 d)$ (with $d \in \mathbb{N}^{*}$) always   by definition of parameter $\alpha$ given in (\ref{Alpha}) we have $\alpha=1/2$. Then our method  breaks down since the integral $\ds{\int_{0}^{t} \frac{ d \tau  }{(t-\tau)^{2 \alpha} \, \tau^{2 \alpha}} }$ diverges. 
\end{Remarque}	

Let us continue with the proof of this theorem. With these estimates on the terms $(a)$ and $(b)$ above, for all $t>0$ and for $\vert x \vert$  large enough we can write 
\begin{eqnarray*}
& & \left\vert  \int_{0}^{t} K_{m,n}(t-\tau) \ast (u^{k}\partial_x u)(\tau, x) d\tau \right\vert   \leq    C_\eta \frac{\mathfrak{C}^{2}_{0}(u)\, e^{c_\eta\, t}}{k+1}\, t^{4 \alpha -1}   \left(  \sup_{0<\tau < t}  \Vert u(\tau, \cdot) \Vert_{H^s}\right)^{k-1}  \, \frac{1}{1+\vert x \vert^{n +2}}\\
&\leq & C_\eta \frac{\mathfrak{C}^{2}_{0}(u)\, e^{c_\eta\, t}}{k+1}\, t^{4 \alpha -1}   \left(  \sup_{0<\tau < t}  \Vert u(\tau, \cdot) \Vert_{H^s}\right)^{k-1}  \, \frac{1}{\vert x \vert^{n +2}} =  \frac{\mathfrak{C}_{1}(u,t)}{\vert x \vert^{n+2}}, 
\end{eqnarray*}
hence, as $\gamma \leq n+2$ we get (\ref{eq19}). \\ 

Now, for the expression $R_1(t,x)$ given in (\ref{eq14}) we set 
\begin{equation}\label{termR}
R(t,x)= R_1(t,x)+ \int_{0}^{t} K_{m,n}(t-\tau) \ast (u^{k}\partial_x u)(\tau, x) d\tau,    
\end{equation}
and then we write 
\begin{equation}\label{Asimp}
u(t,x)= K_{m,n}(t,x)\left( \int_{\R} u_0 (y)dy\right) + R(t,x), 
\end{equation}
where, by estimates (\ref{eq14}) and  (\ref{eq19})  we have the estimate
\begin{equation}\label{estimR}
   \vert R(t,x)\vert \leq \frac{C_1(u_0,u,t)}{\vert x \vert^{\gamma}}, \quad \vert x \vert \to +\infty.  
\end{equation}
With this information, we are able to prove the points $1)$ and $2)$ stated in Theorem \ref{Th-Decay2}. 
\begin{enumerate}
\item[$1)$] If the initial datum verifies $\ds{\int_{\Rt}u_0(y)dy=0}$, then by (\ref{Asimp}) we get $\ds{\vert u(t,x)\vert = \vert R(t,x)\vert }$. From this identity  and by estimate (\ref{estimR}) (recalling that $\gamma=n+1+\varepsilon$)  we obtain the desired estimate (\ref{Decay1Imp}). 
\item[$2)$]  We get back to the identity (\ref{Asimp}), where we assume now that the initial datum verifies $\ds{\int_{\Rt}u_0(y)dy\neq0}$. From this identity we write: 
\begin{eqnarray*}
	\vert u(t,x) \vert &=& \left\vert K_{m,n}(t,x)\left( \int_{\R} u_0 (y)dy\right) + R(t,x) \right\vert = \left\vert K_{m,n}(t,x)\left( \int_{\R} u_0 (y)dy\right) - (- R(t,x))  \right\vert\\
	&\geq & \vert K_{m,n}(t,x) \vert \left\vert \int_{\R} u_0 (y)dy\ \right\vert- \vert R(t,x)\vert. 
\end{eqnarray*}
At this point, we have the following estimate for the kernel $K_{m,n}(t,x)$. 
\begin{Lemme} For $t>0$ fix, there exists a quantity $M=M(t)>0$ such that for all $\vert x \vert > M$ we  have $\ds{\vert K_{m,n}(t,x) \vert \geq \frac{\mathfrak{c}_{\eta}\, t}{2\vert x \vert^{n+1}}}$, for a constant $\mathfrak{c}_{\eta}>0$ depending on $\eta>0$.  
\end{Lemme}	 
\pv  For $n=1$,  by identity (\ref{eq02}) we write 
$$ \vert K_{m,1}(t,x) \vert = \left\vert - \frac{2 \eta t }{(2 \pi i x)^2} + I_1 \right\vert = \left\vert \frac{2 \eta t }{(2 \pi i x)^2}  - I_1 \right\vert \geq \left\vert \frac{2 \eta t }{(2 \pi i x)^2} \right\vert - \vert I_1 \vert= \frac{\mathfrak{c}_\eta\, t}{\vert x \vert^2}-\vert I_1 \vert.$$
Moreover, by (\ref{I1}) we have $\ds{\vert I_1 \vert = o \left( 1/ \vert x \vert^2 \right)}$ when $\vert x \vert \to +\infty$, and then for the quantity $\mathfrak{c}_\eta \,t / 2>0$ there exists $M>0$ such that for all $\vert x \vert >M$ we have $\ds{\vert I_1 \vert \leq \frac{\mathfrak{c}_\eta\, t }{2 \vert x \vert^2}}$. Thus, getting back to the previous estimate we can write $\ds{\frac{\mathfrak{c}_\eta\, t}{\vert x \vert^2}-\vert I_1 \vert \geq \frac{\mathfrak{c}_\eta\, t }{2 \vert x \vert^2}}$ and we obtain the estimate from below  $\ds{\vert K_{m,1}(t,x) \vert  \geq \frac{\mathfrak{c}_\eta\, t }{2 \vert x \vert^2}}$ for all $\vert x \vert >M$. The case $n=2$  follows the same lines with  estimates  (\ref{Estim2}) and (\ref{I2}). Moreover, the case $n\geq 3$ follows the same argument with estimates  (\ref{Estim3}) and (\ref{I3}).  \finpv 
Once we dispose of this lemma, for all $\vert x \vert >M$ we can write $\ds{ \vert u(t,x) \vert \geq \frac{\mathfrak{c}_{\eta}\, t}{2\vert x \vert^{n+1}} \left\vert \int_{\R} u_0 (y)dy\ \right\vert - \vert R(t,x)\vert.}$ Thereafter,  remark that by estimate (\ref{estimR}), with $\gamma=n+1+\varepsilon$ and  $0<\varepsilon\leq 1$, we have $\ds{\vert R(t,x)\vert = o \left(\frac{1}{\vert x \vert^{n+1} }\right)}$   when $\vert x \vert \to +\infty$.  Then,  for the quantity $\ds{\frac{\mathfrak{c}_\eta\, t}{4} \left\vert \int_{\R} u_0 (y)dy \right\vert>0}$, there exists quantity $N=N(\eta, u_0,t)>0$ such that for all $\vert x \vert >N$ we have  $\ds{\vert R(t,x)\vert \leq \frac{\mathfrak{c}_\eta\, t}{4} \left\vert \int_{\R} u_0 (y)dy \right\vert \frac{1}{\vert x \vert^{n+1}}}$. Thus, for all $\vert x \vert > \max(M,N)$ we finally obtain the estimate from below 
$$ \vert u(t,x) \vert \geq \frac{\mathfrak{c}_\eta\, t }{4\, \vert x \vert^{n+1}} \left\vert \int_{\R} u_0 (y)dy\ \right\vert,$$
hence, setting the constant $\ds{C_2= \mathfrak{c}_\eta\, t / 4 >0}$,    the desired estimate (\ref{DecayOp}) follows.   Theorem \ref{Th-Decay2} is now proven. \finpv
\end{enumerate}

\section{Average decaying properties: proof of Theorem \rd{\ref{Th-Averaged-Decay}}}\label{AVD}
To prove this theorem we will follow some of the ideas of the proof of Theorem \ref{Th-Decay}. We assume that for $1<p<+\infty$ and for $0<\gamma<1$ the initial datum $u_0 \in H^s(\R)$ (with $s>3/2$) verifies $u_0 \in L^{p}_{w_\gamma}(\R)$. First, we will construct a solution $$u_1 \in \mathcal{C}([0,+\infty[,H^s(\R))\cap L^{\infty}_{loc}\left(]0,+\infty[, L^{p}_{w_\gamma}(\R), t^{\alpha}\,dt\right),$$ of equation (\ref{Equation-int-aux}). Here, for the parameter $0<\alpha \leq 1/2$ given in (\ref{Alpha}) the weight in the temporal variable $t^\alpha$ is essentially technical and it will be useful to carry up all our estimates.\\

We start by the local in time existence of the solution $u_1(t,x)$ and for this, for a time $0<T_0<+\infty$ small enough, we will solve the equation (\ref{Equation-int-aux}) in the Banach space 
$$ F_{T_0}= \mathcal{C}([0,T_0), H^s(\R)) \cap L^{\infty}(]0,T_0], L^{p}_{w_\gamma}(\R), t^\alpha\, dt),$$ with the norm 
$$ \Vert u \Vert_{F_T}= \sup_{0\leq t \leq T_0} \Vert u(t,\cdot)\Vert_{H^s}+ \sup_{0<t\leq T_0} t^{\alpha} \Vert u (t,\cdot)\Vert_{L^{p}_{w_\gamma}}.$$
For the first term in the right-hand side of (\ref{Equation-int-aux}), recall that the quantity $\ds{\sup_{0\leq t \leq T} \Vert K_{m,n}(t,\cdot)\ast u_0\Vert_{H^s}}$ was estimated in (\ref{estim-lin1}) and then it remains to estimate the quantity $\ds{\sup_{0<t\leq T_0} t^{\alpha} \Vert K_{m,n}(t,\cdot)\ast u_0\Vert_{L^{p}_{w_\gamma}}}$. This estimate bases on the following technical lemma. 
\begin{Lemme}\label{Lema-tech-Lw} For $\psi \in L^{p}_{w_\gamma}(\R)$, and for all $t>0$ we have $\ds{\Vert K_{m,n}(t,\cdot)\ast \psi \Vert_{L^{p}_{w_\gamma}} \leq C_\eta \frac{e^{c_\eta\, t}}{t^{\alpha}} \Vert \psi\Vert_{L^{p}_{w_\gamma}}}$. 
\end{Lemme}
\pv By estimate  (\ref{Kernel-estimates}), for $t>0$ and $x\in \R$ fix we have the following pointwise estimate 
$$\vert K_{m,n}(t,\cdot)\ast \psi (x) \vert \leq (\vert K_{m,n}(t,\cdot)\vert \ast \vert \varphi \vert) (x) \leq C_\eta \frac{e^{c_\eta\, t}}{t^\alpha} \left[ \left( \frac{1}{1+\vert \cdot \vert^{n+1}}\right) \ast \vert \psi \vert \right] (x).$$
Moreover, since the function $\ds{\frac{1}{1+\vert \cdot \vert^{n+1}}}$ belongs to the space $L^1(\R)$  and it is a radially decreasing function then, by the Hardy–Littlewood maximal function operator $\mathcal{M}$ (see the Section $2$ of the book \cite{GrafakosC} for a definition) we can write the pointwise estimate 
$$ C_\eta \frac{e^{c_\eta\, t}}{t^\alpha} \left[ \left( \frac{1}{1+\vert \cdot \vert^{\beta}}\right) \ast \vert \psi \vert  \right] (x) \leq C_\eta  \frac{e^{c_\eta\, t}}{t^\alpha} \left\Vert \frac{1}{1+\vert \cdot\vert^{n+1}}\right\Vert_{L^1} \mathcal{M}_{\vert \psi \vert}(x) \leq C_\eta \frac{e^{c_\eta\, t}}{t^\alpha} \mathcal{M}_{\vert \psi \vert}(x).$$
Thus we obtain  
$$ \Vert K_{m,n}(t,\cdot)\ast \psi \Vert_{L^{p}_{w_\gamma}} \leq  C_\eta \frac{e^{c_\eta\, t}}{t^\alpha} \left\Vert \mathcal{M}_{\vert \psi \vert} \right\Vert_{L^{p}_{w_\gamma}}.$$ 
Now, by Lemma $1$ of \cite{PFPG} we have  that, for $0<\gamma<1$ and for $1<p<+\infty$, the weight $\ds{w_\gamma(x)=\frac{1}{(1+\vert x \vert)^\gamma}}$ belongs to the Muckenhoupt class $\mathcal{A}_{p}(\R)$ (see the book \cite{Grafakos} for a definition). Moreover, by well-known properties of the Muckenhoupt class \cite{Grafakos} we have that  the Hardy–Littlewood maximal function operator $\mathcal{M}$ is bounded in $L^{p}_{w_\gamma}(\R)$  and finally we can write 
$\ds{ C_\eta \frac{e^{c_\eta\, t}}{t^\alpha} \left\Vert \mathcal{M}_{\vert \psi \vert} \right\Vert_{L^{p}_{w_\gamma}} \leq C_\eta \frac{e^{c_\eta\, t}}{t^\alpha} \Vert  \psi \Vert_{L^{p}_{w_\gamma}}}$. \finpv 
By this lemma, with $\psi = u_0$, and by estimate (\ref{estim-lin1}) we get:
\begin{equation}\label{estim-lin2}
\Vert K_{m,n}(t,\cdot)\ast u_0 \Vert_{F_{T_0}} \leq C_\eta \, e^{c_\eta\, T}(\Vert u_0 \Vert_{H^s}+\Vert u_0 \Vert_{L^{p}_{w_\gamma}}).  
\end{equation}
We study now the second term in the right-hanf side of (\ref{Equation-int-aux}). As before, we know that  the  quantity $\ds{\sup_{0<t\leq T_0} \left \Vert  \int_{0}^{t} K_{m,n}(t-\tau, \cdot)\ast (u^{k}_{1} \partial_x u_1)(\tau, \cdot) d \tau \right\Vert_{H^s}}$ was estimated in (\ref{estim-nonlin1}) so  it remains to estimate the quantity $\ds { \sup_{0<t\leq T_0} t^{\alpha} \left \Vert  \int_{0}^{t} K_{m,n}(t-\tau, \cdot)\ast (u^{k}_{1} \partial_x u_1)(\tau, \cdot) d \tau )\right\Vert_{L^{p}_{w_\gamma}}}$. For $0<t<T_0$ fix, and by Lemma \ref{Lema-tech-Lw}  we write  \vspace{2mm}
\begin{eqnarray*}
& &  t^{\alpha} \left \Vert  \int_{0}^{t} K_{m,n}(t-\tau, \cdot)\ast (u^{k}_{1} \partial_x u_1)(\tau, \cdot) d \tau \right\Vert_{L^{p}_{w_\gamma}} \leq  t^{\alpha}  \int_{0}^{t}  \left \Vert   K_{m,n}(t-\tau, \cdot)\ast (u^{k}_{1} \partial_x u_1)(\tau, \cdot)  )\right\Vert_{L^{p}_{w_\gamma}} d \tau \\
&\leq & C_\eta \,  t^{\alpha}  \int_{0}^{t} \frac{e^{c_\eta (t-\tau)}}{(t-\tau)^{\alpha}} \Vert u^{k}_{1} \partial_x u_1)(\tau, \cdot)  \Vert_{L^{p}_{w_\gamma}} d\tau \\
&\leq &  C_\eta \,   e^{c_\eta T_0}   t^{\alpha}  \int_{0}^{t} \frac{1}{(t-\tau)^{\alpha}} \Vert u_1 (\tau, \cdot)\Vert_{L^{p}_{w_\gamma}} \Vert u^{k-1}_{1}(t,\cdot) \Vert_{L^{\infty}} \Vert \partial_x u_1 (t,\cdot) \Vert_{L^{\infty}} d \tau =  (a). 
 \end{eqnarray*}
But, recalling that as $s-1>1/2$  then the spaces $H^s(\R)$ and $H^{s-1}(\R)$ embed in the space $L^{\infty}(\R)$, then we have 
\begin{eqnarray}\label{estim-gron} \nonumber
(a)&\leq & \leq  C_\eta\,  e^{c_\eta\, T_0}   t^{\alpha}  \int_{0}^{t} \frac{1}{(t-\tau)^{\alpha}} \Vert u_1 (\tau, \cdot)\Vert_{L^{p}_{w_\gamma}} \Vert u_1(\tau, \cdot)\Vert^{k}_{H^s} d \tau \\
&\leq &  C_\eta\, e^{c_\eta\, T_0}   t^{\alpha}  \int_{0}^{t} \frac{1}{(t-\tau)^{\alpha} \tau^{\alpha}} (\tau^{\alpha}\, \Vert u_1 (\tau, \cdot)\Vert_{L^{p}_{w_\gamma}}) \Vert u_1(\tau, \cdot)\Vert^{k}_{H^s} d \tau =(b).
\end{eqnarray} 
Moreover, recalling the definition of the norm $\Vert \cdot \Vert_{F_T}$ given above,  we get 
\begin{equation*}
(b)  \leq    C_\eta \,   e^{c_\eta \, T_0}   t^{\alpha} \left( \int_{0}^{t} \frac{1}{(t-\tau)^{\alpha} \tau^{\alpha}} d\tau \right) \Vert u_1 \Vert^{k+1}_{F_T} 
\leq  C_\eta\,  e^{c_\eta \,  T_0}   t^{1-\alpha}  \Vert u_1 \Vert^{k+1}_{F_T}
\leq C_\eta\,  e^{c_\eta \, T_0} \,  T^{1-\alpha}_{0}  \Vert u_1 \Vert^{k+1}_{F_T}.     
\end{equation*}
With this estimate and by estimate (\ref{estim-nonlin1}) we finally write 
\begin{equation}\label{estim-nonlin2}
\left\Vert \int_{0}^{t} K_{m,n}(t-\tau, \cdot)\ast (u^{k}_{1} \partial_x u_1)(\tau, \cdot) d \tau \right\Vert_{F_{T_0}} \leq C_\eta\, e^{c_\eta\, T_0} \,  T^{1-\alpha}_{0}  \Vert u_1 \Vert^{k+1}_{F_T}
\end{equation}
Thus, by estimates (\ref{estim-lin2}) and (\ref{estim-nonlin2}) the existence and uniqueness of  a solution $u_1 \in F_{T_0}$ for  the equation (\ref{Equation-int-aux}) follow from well-know arguments. \\

Recall that by Proposition \ref{Prop-global} we have $u_1 \in \mathcal{C}([0,+\infty[,H^s(\R))$, and the it remains to verify that the quantity  $g(t)= t^{\alpha}\Vert u_1(t,\cdot)\Vert_{L^{p}_{w_\gamma}}$  does not explode in a finite time. For this, for a time $0<T<+\infty$ large enough and for $0<t<T$,  by Lemma \ref{Lema-tech-Lw} and by estimate (\ref{estim-gron}) we can write 
\begin{eqnarray*}
h(t)&\leq & t^{\alpha} \Vert K_{m,n}(t,\cdot)\ast u_0 \Vert_{L^{p}_{w_\gamma}}+t^{\alpha}\left\Vert \int_{0}^{t} K_{m,n}(t-\tau, \cdot)\ast (u^{k}_{1} \partial_x u_1)(\tau, \cdot) d \tau \right\Vert_{L^{p}_{w_\gamma}}\\
&\leq & C_\eta \,  e^{c_\eta\,  T} \Vert u_0 \Vert_{L^{p}_{w_\gamma}}+   C_\eta \,  e^{c_\eta \, T}   t^{\alpha}  \int_{0}^{t} \frac{1}{(t-\tau)^{\alpha} \tau^{\alpha}} (\tau^{\alpha}\, \Vert u_1 (\tau, \cdot)\Vert_{L^{p}_{w_\gamma}}) \Vert u_1(\tau, \cdot)\Vert^{k}_{H^s} d \tau\\
&\leq & C_\eta\, e^{c_\eta\, T} \Vert u_0 \Vert_{L^{p}_{w_\gamma}}+ C_\eta \,  e^{c_\eta \, T}   \left(\sup_{0\leq \tau \leq T} \Vert u_1(\tau,\cdot)\Vert^{k}_{H^s}\right) t^{\alpha}\,  \int_{0}^{t} \frac{1}{(t-\tau)^{\alpha} \tau^{\alpha}} g(\tau) d \tau\\
&\leq& C_1(u_0, T)+C_2(u,T) \, t^{\alpha}\,\int_{0}^{t} \frac{1}{(t-\tau)^{\alpha} \tau^{\alpha}} g(\tau) d \tau,
\end{eqnarray*} and we conclude following the same arguments done at the end of the proof of Theorem \ref{Th-Decay} where  we have treated the cases $0<\alpha <1/2$ (using the Lemma \ref{Lemme-growall-tech})  and $\alpha =1/2$ (using the Lemma \ref{Lemme-growall-tech-2}) separately. \\

We have thus $u_1 \in \mathcal{C}([0,+\infty[,H^s(\R))\cap L^{\infty}_{loc}\left(]0,+\infty[, L^{p}_{w_\gamma}(\R), t^{\alpha}\,dt\right)$. Then, always  by uniqueness of solution $u(t,x)$ in the space $\mathcal{C}([0,+\infty[,H^s(\R))$ we have the identity $u=u_1$ and then, the solution $u(t,x)$ verifies $\ds{u \in L^{\infty}_{loc}\left(]0,+\infty[, L^{p}_{w_\gamma}(\R), t^{\alpha}\,dt\right)}$. This theorem in proven. \finpv
\section{Pointwise growing properties: proof of Theorem \ref{Th-Grow}}\label{PWG}
As the initial datum $u_0 \in \dot{H}^{1}(\R)$ verifies $\vert u_0 (x) \vert \leq C_0 (1+\vert x \vert)^{\gamma}$, then  we get 
$\ds{\left\Vert \frac{1}{(1+\vert \cdot \vert)^\gamma}u_0 \right\Vert_{L^{\infty}}\leq C_0<+\infty}$. Thus, we will solve the integral problem 
\begin{equation}\label{Equation-Int-2}
u(t,x)= K_{m,n}(t,\cdot)\ast u_0(x)+ \int_{0}^{t} K_{m,n}(t-\tau) \ast (u \, \partial_x u)(\tau, x) d\tau,
\end{equation} 
in the Banach space $\ds{G_T=\{ u \in \dot{H}^{1}(\R): \Vert u \Vert_{G_T}<+\infty\}}$, where, for a time $0<T<+\infty$ arbitrary large and fix, and moreover, for the parameter $0<\alpha \leq 1/2$ given in (\ref{Alpha}), we define the norm 
$$ \Vert u \Vert_{G_T}= \sup_{0<t\leq T} \Vert u(t,\cdot)\Vert_{\dot{H}^{1}}+ \sup_{0<t\leq T} t^{\alpha}\left\Vert \frac{1}{(1+\vert \cdot \vert)^\gamma}\,  u(t,\cdot)\right\Vert_{L^{\infty}}.$$ 
Remark that the first term of this norm is technical and it will be useful to treat the second term in the right-hand side of the integral formulation above. On the other hand, observe that the second term in this norm characterizes the pointwise spatial growing of the solution $u(t,x)$. Finally, the weight in time $t^{\alpha}$ is always a technical requirement to carry up of computations.  \\

We start by estimating the first term in the right-hand side of (\ref{Equation-Int-2}). Recall  first that  for all $t>0$ and all $\xi \in \R$  by estimates (\ref{FK2}) and (\ref{FK3}) we have $\ds{\vert \widehat{K_{m,n}}(t,\xi) \vert \leq c}$, and then we write 
\begin{equation}\label{estimg1}
\Vert K_{m,n}(t,\cdot)\ast u_0 \Vert_{\dot{H}^{1}} \leq c\,  \Vert u_0 \Vert_{\dot{H}^{1}}.    
\end{equation}
Next,  by the kernel estimates (\ref{Kernel-estimates}), for  $0<t\leq T$ and $x\in \R$ fix we have the pointwise estimates:
\begin{eqnarray*}
& & \vert K_{m,n}(t,\cdot)\ast u_0 (x)\vert \leq \int_{\R} \vert K_{m,n}(t,x-y)\vert \vert u_0 (y)\vert dy \leq C_\eta \frac{e^{c_\eta \, t}}{t^{\alpha}} \int_{\R} \frac{\vert u_0(y)\vert}{1+\vert x-y \vert^{n+1}} dy \\
&\leq & C_\eta \frac{e^{c_\eta \, t}}{t^{\alpha}} \int_{\R} \frac{(1+\vert y \vert)^{\gamma}\vert u_0(y)\vert}{(1+\vert y \vert)^{\gamma}(1+\vert x-y \vert^{n+1})} dy \leq  C_\eta \frac{e^{c_\eta \, t}}{t^{\alpha}}\, C_0  \int_{\R} \frac{(1+\vert y \vert)^{\gamma}}{1+\vert x-y \vert^{n+1}} dy, 
\end{eqnarray*} where the last expression computes down as $\ds{\int_{\R} \frac{(1+\vert y \vert)^{\gamma}}{1+\vert x-y \vert^{n+1}} dy \leq c_\gamma (1+\vert x \vert)^{\gamma}}$. Indeed, we write 
\begin{eqnarray*}
\int_{\R} \frac{(1+\vert y \vert)^{\gamma}}{1+\vert x-y \vert^{n+1}} dy &\leq &c_\gamma \int_{\R} \frac{(1+\vert x-y \vert)^{\gamma}+(1+\vert x \vert)^\gamma}{1+\vert x-y \vert^{n+1}} dy \\
&\leq & \leq c_\gamma \int_{\R} \frac{(1+\vert x-y \vert)^{\gamma}}{1+\vert x-y \vert^{n+1}} dy + c_\gamma (1+\vert x \vert)^{\gamma} \int_{\R} \frac{1}{1+\vert x-y \vert^{n+1}} dy \\
&\leq & c_\gamma (1+\vert x \vert)^\gamma \left( \int_{\R}  \frac{(1+\vert x-y \vert)^{\gamma}}{1+\vert x-y \vert^{n+1}} dy + \int_{\R} \frac{1}{1+\vert x-y \vert^{n+1}} dy \right), 
\end{eqnarray*}
where, as $n\geq 1$, and moreover, as we have $0<\gamma<1/2$, then  the both integrals above converge.  Thus, getting back to the previous estimate we get 
\begin{equation}\label{estimg2}
t^{\alpha} \left\Vert \frac{1}{(1+\vert \cdot \vert)^\gamma} K_{m,n}(t,\cdot)\ast u_0 \right\Vert_{L^{\infty}} \leq C_\eta \,e^{c_\eta \,T}\, C_0.    
\end{equation}
Thereafter, by (\ref{estimg1}) and (\ref{estimg2}) we  obtain 
\begin{equation}\label{estimglin}
\Vert K_{m,n}(t,\cdot)\ast u_0 \Vert_{G_T}\leq C_{\eta,\gamma}\, e^{c_\eta \,  T} (\Vert u_0 \Vert_{\dot{H}^{1}} + C_0). 
\end{equation}
We study now the second term in the right-hand side of (\ref{Equation-Int-2}). For the first term in the norm $\Vert \cdot \Vert_{G_T}$ we   have 
\begin{eqnarray*}
 \left\Vert \int_{0}^{t} K_{m,n}(t-\tau, \cdot)\ast (u \partial_x u)(\tau, \cdot) d \tau \right\Vert_{\dot{H}^{1}}& =& \frac{1}{2} \left\Vert \int_{0}^{t} K_{m,n}(t-\tau, \cdot)\ast \partial_x  (u^2)(\tau, \cdot) d \tau \right\Vert_{\dot{H}^{1}}\\
&\leq & \int_{0}^{t} \Vert K_{m,n}(t-\tau, \cdot) \ast u^2(\tau, \cdot) \Vert_{\dot{H}^{2}} d \tau. 
\end{eqnarray*}
Then, in point $2)$ of Lemma \ref{Lemma-tech-1} we set the parameters $s_1=3/2$, $s_2=1/2$ and moreover $\psi = u^2$, hence we can write 
\begin{eqnarray*}
 \int_{0}^{t} \Vert K_{m,n}(t-\tau, \cdot) \ast u^2(\tau, \cdot) \Vert_{\dot{H}^{2}} d \tau &\leq & C^{'}_{\eta}  \int_{0}^{t} \frac{e^{c^{'}_\eta (t-\tau)}}{(t-\tau)^{\alpha / 2}} \Vert u^2(\tau, \cdot)\Vert_{\dot{H}^{3/2}} d \tau \\
 &\leq& C^{'}_{\eta}\, e^{c^{'}_{\eta} \,  T} \int_{0}^{t} \frac{1}{(t-\tau)^{\alpha / 2}} \Vert u^2(\tau, \cdot)\Vert_{\dot{H}^{3/2}} d \tau.
\end{eqnarray*}
Moreover, by the  product laws of the homogeneous Sobolev spaces  we have $\ds{\Vert u^2(\tau, \cdot) \Vert_{\dot{H}^{3/2}}\leq c \Vert u(\tau, \cdot) \Vert^{2}_{\dot{H}^{1}}}$. With this estimate, and recalling the definition of the norm $\Vert \cdot\Vert_{G_T}$, from the last expression we get 
\begin{eqnarray*}
C^{'}_{\eta}\, e^{c^{'}_{\eta} \, T} \int_{0}^{t} \frac{1}{(t-\tau)^{\alpha / 2}} \Vert u^2(\tau, \cdot)\Vert_{\dot{H}^{3/2}} d \tau & \leq & C^{'}_{\eta}\, e^{c^{'}_{\eta} \, T} \int_{0}^{t} \frac{1}{(t-\tau)^{\alpha / 2}}  \Vert u(\tau, \cdot)\Vert^{2}_{\dot{H}^{1}} d \tau \\
&\leq & C^{'}_{\eta}\, e^{c^{'}_{\eta} \, T}\left( \sup_{0\leq \tau \leq T} \Vert u(\tau, \cdot)\Vert_{\dot{H}^{1}}\right)^2 \, \int_{0}^{t}\frac{d \tau}{((t-\tau)^{\alpha / 2})} \\
&\leq &C^{'}_{\eta}\, e^{c^{'}_{\eta} \, T}T^{1-\alpha/2} \Vert u \Vert^{2}_{G_T}. 
\end{eqnarray*}
By these estimates we have 
\begin{equation}\label{estimg-nonlin1}
\sup_{0\leq t \leq T} \left\Vert \int_{0}^{t} K_{m,n}(t-\tau, \cdot)\ast (u \partial_x u)(\tau, \cdot) d \tau \right\Vert_{\dot{H}^{1}} \leq C^{'}_{\eta}\, e^{c^{'}_{\eta} \, T} T^{1-\alpha/2} \Vert u \Vert^{2}_{G_T}.
\end{equation}
We estimate now the secod term in the norm $\Vert \cdot \Vert_{G_T}$. By the kernel estimates (\ref{Kernel-estimates}), and moreover, by the definition of the norm $\Vert \cdot \Vert_{G_T}$, for $0<t\leq T$ and $x\in \R$ fix we have the following pointwise estimates: 
\begin{eqnarray*}
& & \left\vert \int_{0}^{t}  K_{m,n}(t-\tau, \cdot) \ast (u \partial_x u)(\tau, x) d \tau \right\vert \leq  \int_{0}^{t} \int_{\R} \vert  K_{m,n}(t-\tau,x-y)\vert \vert  u(\tau,y) \vert \vert  \partial_y u (\tau, y)\vert  dy\, d \tau\\
&\leq &C_\eta  \int_{0}^{t} \frac{e^{c_\eta (t-\tau)}}{(t-\tau)^{\alpha}}\, \left(  \int_{\R} \frac{1}{1+\vert x-y\vert^{n+1}} \vert  u(\tau,y) \vert \vert  \partial_y u (\tau, y)\vert  dy\right)\, d \tau\\
&\leq & C_\eta e^{c_\eta \,  T} \int_{0}^{t} \frac{1}{(t-\tau)^{\alpha}}\,\left(\int_{\R} \frac{(1+\vert y\vert)^\gamma}{1+\vert x-y\vert^{n+1}} \frac{\vert  u(\tau,y) \vert}{(1+\vert y \vert)^{\gamma}} \vert  \partial_y u (\tau, y)\vert  dy \right)\, d \tau\\
&\leq & C_\eta e^{c_\eta \, T} \int_{0}^{t} \frac{1}{(t-\tau)^{\alpha}} \left\Vert \frac{1}{(1+\vert \cdot\vert)^\gamma} u(\tau, \cdot)\right\Vert_{L^{\infty}} \left(\int_{\R} \frac{(1+\vert y\vert)^\gamma}{1+\vert x-y\vert^{n+1}}  \vert  \partial_y u (\tau, y)\vert  dy \right)\, d \tau\\
&\leq & C_\eta e^{c_\eta \, T} \int_{0}^{t} \frac{1}{(t-\tau)^{\alpha}\,\tau^{\alpha}} \left(\tau^{\alpha}\left\Vert \frac{1}{(1+\vert \cdot\vert)^\gamma} u(\tau, \cdot)\right\Vert_{L^{\infty}}\right) \left(\int_{\R} \frac{(1+\vert y\vert)^\gamma}{1+\vert x-y\vert^{n+1}}  \vert  \partial_y u (\tau, y)\vert  dy \right)\, d \tau\\
 &\leq & C_\eta \, e^{c_\eta \, T}\left(\sup_{0\leq \tau \leq T} \tau^{\alpha}\left\Vert \frac{1}{(1+\vert \cdot\vert)^\gamma} u(\tau, \cdot)\right\Vert_{L^{\infty}}\right) \int_{0}^{t} \frac{1}{(t-\tau)^{\alpha}\,\tau^{\alpha}} \left(\int_{\R} \frac{(1+\vert y\vert)^\gamma}{1+\vert x-y\vert^{n+1}}  \vert  \partial_y u (\tau, y)\vert  dy \right)\, d \tau\\
 &\leq & C_\eta \, e^{c_\eta \, T}\, \Vert u \Vert_{G_T} \int_{0}^{t} \frac{1}{(t-\tau)^{\alpha}\,\tau^{\alpha}} \left(\int_{\R} \frac{(1+\vert y\vert)^\gamma}{1+\vert x-y\vert^{n+1}}  \vert  \partial_y u (\tau, y)\vert  dy \right)\, d \tau =(a).
\end{eqnarray*} At this point, we must study the integral in the spatial variable. Applying first  the Cauchy-Schwarz inequalities  we write 
\begin{eqnarray*}
 \int_{\R} \frac{(1+\vert y\vert)^\gamma}{1+\vert x-y\vert^{n+1}}  \vert  \partial_y u (\tau, y)\vert  dy &\leq &\left( \int_{\R} \frac{(1+\vert y\vert)^{2\gamma}}{1+\vert x-y\vert^{2(n+1)}} dy \right)^{1/2} \left(\int_{\R} \vert  \partial_y u (\tau, y)\vert  dy \right)^{1/2} \\
 &\leq & \left( \int_{\R} \frac{(1+\vert y\vert)^{2\gamma}}{1+\vert x-y\vert^{2(n+1)}} dy \right)^{1/2}  \Vert u(\tau, \cdot) \Vert_{\dot{H}^{1}} \leq  c_\gamma (1+\vert x \vert)^\gamma \, \Vert u(\tau, \cdot) \Vert_{\dot{H}^{1}}. 
\end{eqnarray*}
Thus, getting back to the term  $(a)$, by this estimate and always  by the definition of $\Vert \cdot\Vert_{G_T}$,  we obtain
\begin{eqnarray*}
(a)&\leq & C_{\eta,\gamma} (1+\vert x \vert)^{\gamma}\, e^{c_\eta\, T}\, \Vert u \Vert_{G_T} \int_{0}^{t} \frac{1}{(t-\tau)^{\alpha}\,\tau^{\alpha}} \Vert u(\tau, \cdot)\Vert_{\dot{H}^{1}}  d \tau \\
&\leq & C_{\eta,\gamma} (1+\vert x \vert)^{\gamma}\, e^{c_\eta \,T}\, \Vert u \Vert_{G_T} \left( \sup_{0 \leq \tau \leq T} \Vert u(\tau, \cdot)\Vert_{\dot{H}^{1}}\right) \int_{0}^{t} \frac{d \tau}{(t-\tau)^{\alpha}\,\tau^{\alpha}}\\
&\leq &  C_{\eta,\gamma} (1+\vert x \vert)^{\gamma}\, e^{c_\eta\, T}\, \Vert u \Vert^{2}_{G_T}\, t^{1-2\alpha}. 
\end{eqnarray*}
By these estimates we  get 
\begin{equation}\label{estimg-nonlin2}
\sup_{0\leq t \leq T} t^\alpha\left\Vert \left(\frac{1}{(1+\vert \cdot\vert)^\gamma}\right)  \int_{0}^{t}  K_{m,n}(t-\tau, \cdot) \ast (u \partial_x u)(\tau, x) d \tau \right\Vert_{L^{\infty}} \leq  C_{\eta,\gamma}  e^{c_\eta \,T}\, T^{1-\alpha}\, \Vert u \Vert^{2}_{G_T}.
\end{equation}
Finally, by estimates (\ref{estimg-nonlin1}) and (\ref{estimg-nonlin2}) we have 
\begin{equation}\label{estimg-nonlin}
\left\Vert \int_{0}^{t}  K_{m,n}(t-\tau, \cdot) \ast (u \partial_x u)(\tau, x) d \tau \right\Vert_{G_T} \leq C_{\eta,\gamma} \, e^{c_\eta \,  T}\,\max(T^{1-\alpha / 2}, T^{1-\alpha}) \, \Vert u \Vert^{2}_{G_T}.  
\end{equation}
Now, with the estimates (\ref{estimglin}) and (\ref{estim-nonlin2}) we set the quantity $\delta$ as $\ds{\delta = \frac{1}{4 C_{\eta,\gamma} \, e^{c_\eta\,  T}}>0}$, and if  the initial datum verifies $\ds{\Vert u_0 \Vert_{\dot{H}^{1}} + C_0 < \delta}$ then the existence and uniqueness of a solution $u \in G_T$ of  equation (\ref{Equation-int-aux}) follow from standard arguments. Theorem \ref{Th-Grow} is now proven. \finpv  

			 \quad\\[5mm]
		 \begin{flushright}
	 \begin{minipage}[r]{130mm}
	 \textbf{{\large Manuel Fernando Cortez}} (manuel.cortez@epn.edu.ec) \\ 
	 \\
	 Departamento de Matem\'aticas, Escuela Politécnica Nacional,  Ladron de Guevera E11-253n Quito, Ecuador.\\ 
	 \\
	 \\
	 \\
	 \textbf{\large {Oscar Jarr\'in}} (corresponding author: oscar.jarrin@udla.edu.ec)\\
	 \\
	Dirección General  de Investigación  (DGI),
	Universidad de las Américas,
	Calle José Queri s/n y Av. Granados. Bloque 7, Tercer Piso, Quito, Ecuador.
	 \end{minipage}
	 \end{flushright}
\end{document}